\documentstyle{article}
\begin{document}
\font\germ=eufm10
\def\ssl{\hbox{\germ sl}}
\def\slh{\widehat{\ssl_2}}
\makeatletter
\def\aaa{@}
\centerline{\Large\bf Polyhedral Realizations of Crystal Bases for}
\centerline{\Large\bf Integrable Highest Weight Modules}
\vskip12pt
\centerline{Toshiki NAKASHIMA}
\vskip8pt
\centerline{Department of Mathematics,}
\centerline{Sophia University, Tokyo 102-8554, JAPAN}
\centerline{toshiki@mm.sophia.ac.jp}
\makeatother

\renewcommand{\labelenumi}{(\roman{enumi})}
\font\germ=eufm10

\def\al{\alpha}
\def\beneme{\begin{enumerate}}
\def\beq{\begin{equation}}
\def\beqn{\begin{eqnarray}}
\def\beqnn{\begin{eqnarray*}}
\def\bigsl{{\hbox{\fontD \char'54}}}
\def\cd{\cdots}
\def\del{\delta}
\def\Del{\Delta}
\def\ei{e_i}
\def\eit{\tilde{e}_i}
\def\eneme{\end{enumerate}}
\def\ep{\epsilon}
\def\eeq{\end{equation}}
\def\eeqn{\end{eqnarray}}
\def\eeqnn{\end{eqnarray*}}
\def\fit{\tilde{f}_i}
\def\ft{\tilde{f}}
\def\ge{\hbox{\germ g}}
\def\gl{\hbox{\germ gl}}
\def\hom{{\hbox{Hom}}}
\def\ify{\infty}
\def\io{\iota}
\def\kp{k^{(+)}}
\def\km{k^{(-)}}
\def\llra{\relbar\joinrel\relbar\joinrel\relbar\joinrel\rightarrow}
\def\lan{\langle}
\def\lar{\longrightarrow}
\def\lm{\lambda}
\def\Lm{\Lambda}
\def\mapright#1{\smash{\mathop{\longrightarrow}\limits^{#1}}}
\def\nd{\noindent}
\def\nn{\nonumber}
\def\ot{\otimes}
\def\op{\oplus}
\def\opi{\ovl\pi_{\lm}}
\def\ovl{\overline}
\def\plm{\Psi^{(\lm)}_{\io}}
\def\qq{\qquad}
\def\q{\quad}
\def\qed{\hfill\framebox[3mm]{}}
\def\QQ{\hbox{\bf Q}}
\def\qi{q_i}
\def\qii{q_i^{-1}}
\def\ran{\rangle}
\def\rlm{r_{\lm}}
\def\ssl{\hbox{\germ sl}}
\def\slh{\widehat{\ssl_2}}
\def\ti{t_i}
\def\tii{t_i^{-1}}
\def\til{\tilde}
\def\tt{{\hbox{\germ{t}}}}
\def\ttt{\hbox{\germ t}}
\def\uq{U_q(\ge)}
\def\uqm{U^-_q(\ge)}
\def\uqmq{{U^-_q(\ge)}_{\bf Q}}
\def\uqpm{U^{\pm}_q(\ge)}
\def\uqq{U_{\bf Q}^-(\ge)}
\def\uqz{U^-_{\bf Z}(\ge)}
\def\util{\widetilde\uq}
\def\vep{\varepsilon}
\def\vp{\varphi}
\def\vpi{\varphi^{-1}}
\def\xii{\xi^{(i)}}
\def\Xiioi{\Xi_{\io}^{(i)}}
\def\wtil{\widetilde}
\def\what{\widehat}
\def\wpi{\widehat\pi_{\lm}}
\def\ZZ{\hbox{\bf Z}}

\renewcommand{\thesection}{\arabic{section}}
\section{Introduction}
\setcounter{equation}{0}
\renewcommand{\theequation}{\thesection.\arabic{equation}}

Since Kashiwara introduced the theory of crystal base (\cite{K0}) in
1990,
one of the most fundamental problems has been to 
describe the crystal base associated with the given 
integrable highest weight module
as explicitly as possible.
In order to answer this, many kinds of new 
combinatorial objects have been invented, e.g.,
in \cite{KN} some analogues of Young tableaux were introduced
in order to describe the crystal base for classicla Lie algebras and
it is applied to generalize so-called the 
Littlewood-Richardson rule in\cite{N}.
In \cite{KMN1}\cite{KMN2} we gave the new object 'perfect crystals' and 
applying it to describe the crystal bases of affine types, 
moreover, to solve the problems in mathmatical physics.
In \cite{L1}\cite{L2}
 Littelman realized the crystal base for symmetrizable Kac-Moody
Lie algebras by using 'path' and in \cite{NZ} we also have done it 
for the nilpotent subalgebra $\uqm$ 
by the 'polyhedral realization'.

The present paper is devoted to give the explicit feature 
of crystal bases for integrable highest weight modules
in terms of `polyhedral realization'.
Here we introduce the formulation and background of \cite{NZ}
and this paper. 
Let $\uq$ be the quantum algebra associated with 
the Kac-Moody Lie algebra $\ge$ and $\uqm$ be 
the nilpotent subalgebra given by the usual triangular decomposition
as in 2.1 below.
Furthermore, let $(L(\ify),B(\ify))$ be the crystal base
of $\uqm$ (see \cite{K1}).
In \cite{K3}, Kashiwara introduced the remarkable embedding 
of crystals $\Psi_{\io}:B(\ify)\hookrightarrow \ZZ^{\ify}$
where $\io$ is some infinite sequence from the index set $I$ and 
$\ZZ^{\ify}$ is the $\ZZ$-lattice of infinite rank (see 2.5).
In \cite{NZ}, we tried to describe the exact image of the embedding
$\Psi_{\io}$ in $\ZZ^{\ify}$. This can be carried out by the 
unified method, called  the {\it polyhedral realization}.
Due to this method, (under some condition)
we succeeded to present the explicit
form of $B(\ify)$ as the set of lattice points
of some polyhedral convex cone in 
the infinite vector space $\QQ^{\ify}$,
which is defined by the system of inequalities.
This system of inequalities are determined only by the sequence
$\io$ and the Cartan data of $\ge$.
In the present paper, we shall try to give 
the similar decription of the crystal $B(\lm)$, where 
$B(\lm)$ is the crystal (base) of the irreducible 
integrable highest weight module 
$V(\lm)$ with the highest weight $\lm$.

To mention our problem more precisely, let us introduce
the object $R_{\lm}$ which is a crystal consisting of 
the one element $r_{\lm}$ ($\lm$ is a weitght)
(see Example \ref{Example:crystal} (ii) below and also \cite{J}) and 
has the following remarkable property:
The crystal $B(\ify)$ is connnected as a crystal graph 
(see Definition\ref{c-gra}), 
but in general, not so the crystal $B(\ify)\ot R_{\lm}$ is.
Furthermore, the connected component including 
$u_{\ify}\ot r_{\lm}$ is isomorphic to the crystal 
$B(\lm)$, where 
$u_{\ify}$ is the highest weight vector in $B(\ify)$
(Theorem \ref{ify-lm}).
These properties guarantee the existence of the embedding of crystals 
$\Omega_{\lm}:B(\lm)\hookrightarrow B(\ify)\ot R_{\lm}$.
(see also \cite{J}).
Combining $\Omega_{\lm}$ and $\Psi_{\io}$,
we obtain the embedding of crystals
$\Psi^{(\lm)}_{\io}:=(\Psi_{\io}\ot{\rm id})\circ \Omega_{\lm}
:B(\lm)\hookrightarrow \ZZ^{\ify}\ot R_{\lm}$.
Here note that since $R_{\lm}$ consists of one element, as a set 
$\ZZ^{\ify}\ot R_{\lm}$ can be identified with the infinite 
$\ZZ$-lattice $\ZZ^{\ify}$.
Our goal is to give the explicit form of 
${\rm Im}(\Psi^{(\lm)}_{\io})(\cong B(\lm))$
in the infinite $\ZZ$-lattice.
To complete this, we shall introduce the set linear functions
$\Xi_{\io}[\lm]$ (see (\ref{Xi})) which is uniquely determined 
by the Cartan data of $\ge$, the sequence $\io$ and the 
highest weight $\lm$. The set $\Sigma_{\io}[\lm]$ is the 
set of lattice points in the convex polyhedron 
defined by the system of inequalities:
$\vp(\vec x)\geq 0$ ($\vec x\in \ZZ^{\ify}$)
for any $\vp\in \Xi_{\io}[\lm]$.
Finally, we can show ${\rm Im}(\Psi^{(\lm)}_{\io})=\Sigma_{\io}[\lm]$
under the assumption $\Sigma_{\io}[\lm]\ni\vec 0:=(\cd, 0,0,0)$.
We shall apply this to several explicit cases, namely,
arbitrary rank 2 Kac-Moody algebras, $A_n$-case 
and $A_{n-1}^{(1)}$-case. 

Now let us see the organization of this paper. 
The section 2 is devoted to preliminaries and reviews on 
the theory of  crystals and crystal bases. 
In particular, in 2.1 the crystal $R_{\lm}$ will be 
introduced and in 2.4 we shall give the results of \cite{NZ}.
In section 3, we shall introduce the surjective morphism
$\Phi_{\lm}:B(\ify)\ot R_{\lm}\longrightarrow B(\lm)$,
the embeddings $\Omega_{\lm}:B(\lm)\hookrightarrow B(\ify)\ot R_{\lm}$ and 
$\Psi^{(\lm)}_{\io}:B(\lm)\hookrightarrow \ZZ^{\ify}\ot R_{\lm}$.
The section 4 is the main part of this paper. 
In 4.2 we shall construct the polyhedral realization of $B(\lm)$
explicitly and in 4.3 it is applied to give the explicit 
description of the value '$\vep^*_i$'.
In section 5, we treat rank 2 Kac-Moody algebras.
The explicit form of polyhedral realization
of $B(\lm)$ is given by using 'Chebyshev polynomials'.
In section 6, we consider the case $\ge=A_n$.
Furthermore, in this section we shall give an example 
which does not satisfy the `positivity assumption'.
Thus, our conjectural perspective in \cite[3.3]{NZ}, that 
the positivity assumption is satisfied automatically,  
turns out to be invalid.
In section 7, we treat the higher rank affine case
$\ge=A^{(1)}_{n-1}$.

The author would like to thank A.Zelevinsky for 
valuable advices and suggestions.  
He also would like to acknowledge S.Zelikson for notifying him 
that the crystal $R_{\lm}$ and the embedding $\Omega_{\lm}$ 
have already appeared in \cite{J}.
\section{Crystals and Crystal Bases}
\setcounter{equation}{0}
\renewcommand{\thesection}{\arabic{section}}
\renewcommand{\theequation}{\thesection.\arabic{equation}}

\subsection{Definition of crystals}

Let $\ge$ be
a  symmetrizable Kac-Moody algebra over {\bf Q}
with a Cartan subalgebra
$\ttt$, a weight lattice $P \subset \ttt^*$, the set of simple roots
$\{\al_i: i\in I\} \subset \ttt^*$,
and the set of coroots $\{h_i: i\in I\} \subset \ttt$,
where $I$ is a finite index set.
Let $\lan h,\lm\ran$ be the pairing between $\ttt$ and $\ttt^*$,
and $(\al, \beta)$ be an inner product on
$\ttt^*$ such that $(\al_i,\al_i)\in 2{\bf Z}_{\geq 0}$ and
$\lan h_i,\lm\ran={{2(\al_i,\lm)}\over{(\al_i,\al_i)}}$
for $\lm\in\ttt^*$.
Let $P^*=\{h\in \ttt: \lan h,P\ran\subset\ZZ\}$ and
$P_+:=\{\lm\in P:\lan h_i,\lm\ran\in\ZZ_{\geq 0}\}$.
We call an element in $P_+$ a {\it dominant integral weight}.
The quantum algebra $\uq$
is an associative
$\QQ(q)$-algebra generated by the $e_i$, $f_i \,\, (i\in I)$,
and $q^h \,\, (h\in P^*)$
satisfying the usual relations (see e.g.,\cite{K1} or \cite{NZ}).
The algebra $\uqm$ is the subalgebra of $\uq$ generated 
by the $f_i$ $(i\in I)$.

The following definition of a crystal is  the one
slightly modified those in \cite{K3,K4}. But there is no difference 
between their properties.
In what follows we fix a finite index set $I$ and a
weight lattice $P$ as above.

\newtheorem{df2}{Definition}[section]
\begin{df2}
\label{def2.1}
A {\it crystal} $B$ is a set endowed with the following maps:
\begin{eqnarray}
&& wt:B\lar P,\\
&&\vep_i:B\lar\ZZ\sqcup\{-\infty\},\q
  \vp_i:B\lar\ZZ\sqcup\{-\infty\} \q{\hbox{for}}\q i\in I,\\
&&\eit:B\sqcup\{0\}\lar B\sqcup\{0\},
\q\fit:B\sqcup\{0\}\lar B\sqcup\{0\}\q{\hbox{for}}\q i\in I.
\end{eqnarray}
Here $0$ is an
ideal element which is not included in $B$.
These maps must satisfy the following axioms:
for all $b$,$b_1$,$b_2\in B$, we have
\begin{eqnarray}
&&\vp_i(b)=\vep_i(b)+\lan h_i,wt(b)\ran,
\label{vp=vep+wt}\\
&&wt(\eit b)=wt(b)+\al_i{\hbox{ if }}\eit b\in B,
\label{+alpha}\\
&&wt(\fit b)=wt(b)-\al_i{\hbox{ if }}\fit b\in B,\\
&&\eit b_2=b_1 {\hbox{ if and only if }} \fit b_1=b_2,
\label{eeff}\\
&&{\hbox{if }}\vep_i(b)=-\infty,
  {\hbox{ then }}\eit b=\fit b=0,\\
&&\eit(0)=\fit(0)=0.
\end{eqnarray}
\end{df2}

The above axioms allow us to make a crystal $B$ into
a colored oriented graph with the set of colors $I$.

\begin{df2}
\label{c-gra}
The crystal graph of a crystal $B$ is
a colored oriented graph given by
the rule : $b_1\mapright{i} b_2$ if and only if $b_2=\fit b_1$
$(b_1,b_2\in B)$.
\end{df2}

\begin{df2}
\label{df:mor}
\begin{enumerate}
\item
Let $B_1$ and $B_2$ be crystals.
A {\sl strict morphism } of crystals $\psi:B_1\lar B_2$
is a map $\psi:B_1\sqcup\{0\} \lar B_2\sqcup\{0\}$
satisfying the following conditions: $\psi(0)=0$, 
\begin{eqnarray}
&&\hspace{-30pt}wt(\psi(b)) = wt(b),\q \vep_i(\psi(b)) = \vep_i(b),\q
\vp_i(\psi(b)) = \vp_i(b)
\label{wt}\\
&&{\hbox{if }}b\in B_1{\hbox{ and }}\psi(b)\in B_2,\nonumber,
\end{eqnarray}
and the map $\psi: B_1\sqcup\{0\} \lar B_2\sqcup\{0\}$
commutes with all $\eit$ and $\fit$.
\item
An injective strict morphism is called an embedding of crystals.
We call $B_1$ is a subcrystal of $B_2$, if $B_1$ is a subset of $B_2$ and
becomes a crystal itself by restricting the data on it from $B_2$.
\end{enumerate}
\end{df2}

It is well-known that $\uq$ has a Hopf algebra structure.
Then the tensor product of $\uq$-modules has
a $\uq$-module structure.
Consequently, we can consider the tensor product
of crystals:
For crystals $B_1$ and $B_2$, we define their tensor product
$B_1\ot B_2$ as follows:
\begin{eqnarray}
&&B_1\ot B_2=\{b_1\ot b_2: b_1\in B_1 ,\, b_2\in B_2\},\\
&&wt(b_1\ot b_2)=wt(b_1)+wt(b_2),\\
&&\vep_i(b_1\ot b_2)={\hbox{max}}(\vep_i(b_1),
  \vep_i(b_2)-\lan h_i,wt(b_1)\ran),
\label{tensor-vep}\\
&&\vp_i(b_1\ot b_2)={\hbox{max}}(\vp_i(b_2),
  \vp_i(b_1)+\lan h_i,wt(b_2)\ran),
\label{tensor-vp}\\
&&\eit(b_1\ot b_2)=
\left\{
\begin{array}{ll}
\eit b_1\ot b_2 & {\mbox{ if }}\vp_i(b_1)\geq \vep_i(b_2)\\
b_1\ot\eit b_2  & {\mbox{ if }}\vp_i(b_1)< \vep_i(b_2),
\end{array}
\right.
\label{tensor-e}
\\
&&\fit(b_1\ot b_2)=
\left\{
\begin{array}{ll}
\fit b_1\ot b_2 & {\mbox{ if }}\vp_i(b_1)>\vep_i(b_2)\\
b_1\ot\fit b_2  & {\mbox{ if }}\vp_i(b_1)\leq \vep_i(b_2).
\label{tensor-f}
\end{array}
\right.
\end{eqnarray}
Here $b_1\ot b_2$ is just another notation for an ordered pair
$(b_1,b_2)$, and we set $b_1 \ot 0 = 0 \ot b_2=0$.
Note that the tensor product of crystals is
associative, namely, the crystals
$(B_1\ot B_2)\ot B_3$ and $B_1\ot(B_2\ot B_3)$ are isomorphic via
$(b_1\ot b_2)\ot b_3\leftrightarrow b_1\ot (b_2\ot b_3)$.

The example of crystals below will be needed later.
\newtheorem{ex}[df2]{Example}
\begin{ex}
\label{Example:crystal}
\begin{enumerate}
\item
For $i\in I$, the crystal $B_i:=\{(x)_i\,: \, x \in\ZZ\}$ is defined by
\beqnn
&& wt((x)_i)=x \al_i,\qq \vep_i((x)_i)=-x,\qq \vp_i((x)_i)=x,\\
&& \vep_j((x)_i)=-\infty,\qq \vp_j((x)_i)
   =-\infty \q {\rm for }\q j\ne i,\\
&& \til e_j (x)_i=\del_{i,j}(x+1)_i,\qq 
\til f_j(x)_i=\del_{i,j}(x-1)_i,
\eeqnn
\item
(See also \cite{J}) 
Let $R_{\lm}:=\{r_{\lm}\}$ $(\lm\in P)$
be the crystal consisting of one-element given by
\beqnn
&& wt(r_{\lm})=\lm,\q
 \vep_i(r_{\lm})=-\lan h_i,\lm\ran,\q
   \vp_i(r_{\lm})=0,
\q \eit (r_{\lm})=\fit(r_{\lm})=0.
\eeqnn
\end{enumerate}
\end{ex}

\subsection{Crystal $B(\lm)$}

In this subsection we review the crystal $B(\lm)$ for a dominant 
integral weight $\lm\in P_+$,
which is our main object of study.
All the results in this subsection are due to M.Kashiwara \cite{K1}.
Let $V(\lm)$ be the irreducible highest weight module of $\uq$
with the highest weight $\lm\in P_+$.
It can be defined by
\beq
V(\lm)  :=  \uq\left/\sum_i\uq e_i+\sum_i\uq f_i^{\lan h_i,\lm\ran+1}+
\sum_{h\in P^*}\uq(q^h-q^{\lan h,\lm\ran})\right.{.}
\label{vlm}
\eeq
It is well-known that as a $\uqm$-module,
there is the following natural isomorphism:
\beq
V(\lm) \cong  \uqm/\sum_i\uqm f_i^{\lan h_i,\lm\ran+1}.
\label{vlm-uqm}
\eeq
Let $\pi_{\lm}$ be a natural projection $\uqm\longrightarrow V(\lm)$ and
set $u_{\lm}:=\pi_{\lm}(1)$. This is the unique highest weight vector in 
$V(\lm)$ up to constant.

For each $i\in I$, we have the decomposition:
$
V(\lm)=\bigoplus_n f^{(n)}_i({\rm Ker}\,e_i).
$
Using this, we can define the endomorphisms
$\eit$ and $\fit\in {\rm End}(V(\lm))$ by
\beq
\eit(f^{(n)}_iu)=f^{(n-1)}_iu,\,\, {\rm and}\,\,
\fit(f^{(n)}_iu)=f^{(n+1)}_iu\q {\rm for}\,\, u\in\,{\rm Ker }\,\,e_i,
\label{eit-fit}
\eeq
where we understand that $\eit u=0$ for $u\in{\rm Ker}\,e_i$.
Let $A\subset \QQ(q)$ be the subring of rational functions
regular at $q=0$.
We set
\beqn
L(\lm) & :=& \sum_{i_j\in I,l\geq 0}
A\til f_{i_l}\cd \til f_{i_1}u_{\lm},\\
B(\lm) & := &
\{\til f_{i_l}\cd \til f_{i_1}u_{\lm}\,\,
{\rm mod}\,\,qL(\lm)\,|\,i_j\in I,l\geq 0\}
\setminus \{0\}.
\label{def-blm}
\eeqn
The pair $(L(\lm),B(\lm))$ is called {\it crystal base} of $V(\lm)$.
It satisfies the following properties:
\begin{enumerate}
\item
$L(\lm)$ is a free $A$-submodule of $V(\lm)$
and $V(\lm)\cong \QQ(q)\ot_A L(\lm)$.
\item
$B(\lm)$ is a basis of the $\QQ$-vector space $L(\lm)/qL(\lm)$.
\item
$\eit L(\lm)\subset L(\lm)$ and
$\fit L(\lm)\subset L(\lm)$.

By (iii) the $\eit$ and the $\fit$ act on
$L(\lm)/qL(\lm)$ and
\item
$\eit B(\lm)\subset B(\lm)\sqcup\{0\}$ and
$\fit B(\lm)\subset B(\lm)\sqcup\{0\}$.
\item
For $u,v\in B(\lm)$,
$\fit u=v$ if and only if  $\eit v=u$.
\end{enumerate}

We define the weight function $wt:B(\lm)\rightarrow P$
by $wt(b):=\lm-\al_{i_1}-\al_{i_2}-\cd-\al_{i_l}$ for
$b=\til f_{i_l}\cd\til f_{i_1}u_{\lm}$ mod $qL(\lm)\ne0$.
We define integer-valued functions $\vep_i$ and $\vp_i$ on
$B(\lm)$ by
\beq
\vep_i(b):={\rm max} \{k: \eit^k b\ne 0\}, \,\,
\vp_i(b):={\rm max} \{k: \fit^k b\ne 0\}.
\label{vep-vp}
\eeq
It is easy to verify that $B(\lm)$ equipped with the operators
$\eit$ and $\fit$, and with the functions $wt$, $\vep_i$ and $\vp_i$
becomes a crystal.

Let $(L(\ify),B(\ify))$ be the crystal base of the subalgebra
$\uqm$ (see \cite{K1},\cite{NZ}).
Here note that the functions $\vep_i$ and $\vp_i$ are given by 
\beq
\vep_i(b):={\rm max} \{k: \eit^k b\ne 0\}, \,\,
\vp_i(b):=\vep_i(b)+\lan h_i,wt(b)\ran.
\label{vep-vp-ify}
\eeq
It is proved in \cite{K1} that
the natural projection $\pi_{\lm}:\uqm\rightarrow V(\lm)$
sends $L(\ify)$ to $L(\lm)$, and the induced map
$\what\pi_{\lm}:L(\ify)/qL(\ify)\longrightarrow L(\lm)/q L(\lm)$ sends
$B(\ify)$ to $B(\lm)\sqcup\{0\}$.
The map $\wpi$ has the following properties:
\beqn
&&\fit\circ\wpi=\wpi\circ\fit,
\label{fpi=pif}\\
&&\eit\circ\wpi=\wpi\circ\eit,\,\,{\rm if }\,\,\wpi(b)\ne0,
\label{epi=pie}\\
&&\hspace{-10pt}{\hbox{$\wpi:B(\ify)\setminus \{\wpi^{-1}(0)\}\longrightarrow
B(\lm)$ is bijective.}}
\label{wpi-bij}
\eeqn
Although the map $\wpi$ has such nice properties, 
it is not a strict morphism of crystals.
For instance, it does not preserve weights or 
does not necessarily commute with
the action of $\eit$ as in (\ref{epi=pie}).
We shall introduce a new morphism by modifying the map 
$\wpi$ in 3.1.

\subsection{Kashiwara Embedding}

We define a $\QQ(q)$-algebra anti-automorphism $*$ of $\uq$ by:
$q^*=q$, $e^*_i=e_i$,  $f^*_i=f_i$, $(q^h)^*=q^{-h}$.
This antiautomorphism has the properties (see \cite{K3}):
\beq
L(\ify)^*=L(\ify)\,\,{\rm  and}\,\,B(\ify)^*=B(\ify).
\label{*-sta}
\eeq
Then we can define
$\vep^*_i(b):=\vep_i(b^*)$ and $\vp^*_i(b):=\vp_i(b^*)$.

Consider the additive group
\beq
\ZZ^{\ify}
:=\{(\cd,x_k,\cd,x_2,x_1): x_k\in\ZZ
\,\,{\rm and}\,\,x_k=0\,\,{\rm for}\,\,k\gg 0\};
\label{uni-cone}
\eeq
we will denote by $\ZZ^{\ify}_{\geq 0} \subset \ZZ^{\ify}$
the subsemigroup of nonnegative sequences.
To the rest of this section, we fix an infinite sequence of indices
$\io=\cd,i_k,\cd,i_2,i_1$ from $I$ such that
\beq
{\hbox{
$i_k\ne i_{k+1}$ and $\sharp\{k: i_k=i\}=\ify$ for any $i\in I$.}}
\label{seq-con}
\eeq
We can associate to $\io$ a crystal structure
on $\ZZ^{\ify}$ and denote it by $\ZZ^{\ify}_{\io}$ 
(\cite[2.4]{NZ}).

\newtheorem{pro2}[df2]{Proposition}
\begin{pro2}[\cite{K3}, See also \cite{NZ}]
\label{emb}
There is a unique embedding of crystals
\beq
\Psi_{\io}:B(\ify)\hookrightarrow \ZZ^{\ify}_{\geq 0}
\subset \ZZ^{\ify}_{\io},
\label{psi}
\eeq
such that
$\Psi_{\io} (u_{\ify}) = (\cd,0,\cd,0,0)$.
\end{pro2}
We call this the {\it Kashiwara embedding}
which is derived by iterating the following
type of embeddings (\cite{K3}):
\begin{enumerate}
\item
For any $i\in I$, there is a unique embedding of crystals
\beqn
\Psi_i : B(\infty)&\hookrightarrow &B(\infty)\ot B_i,
\label{Psii}
\eeqn
such that $\Psi_i (u_{\ify}) = u_{\ify}\ot (0)_i$.
\item
For any $b\in B(\ify)$, we can write uniquely
 $\Psi_i(b)=b'\ot \fit^m(0)_i$
where $m=\vep^*_i(b)$.
\end{enumerate}

\subsection{Polyhedral Realization of $B(\ify)$}
In this subsection, we recall the main result of \cite{NZ}.

Consider the infinite dimensional vector space
$$
\QQ^{\ify}:=\{\vec{x}=
(\cd,x_k,\cd,x_2,x_1): x_k \in \QQ\,\,{\rm and }\,\,
x_k = 0\,\,{\rm for}\,\, k \gg 0\},
$$
and its dual $(\QQ^{\ify})^*:={\rm Hom}(\QQ^{\ify},\QQ)$.
We will write a linear form $\vp \in (\QQ^{\ify})^*$ as
$\vp(\vec{x})=\sum_{k \geq 1} \vp_k x_k$ ($\vp_j\in \QQ$).

Fix $\io=(i_k)$ as in the previous section.
For $\io$ we set $\kp:={\rm min}\{l:l>k\,\,{\rm and }\,\,i_k=i_l\}$ and
$\km:={\rm max}\{l:l<k\,\,{\rm and }\,\,i_k=i_l\}$ if it exists,
or $\km=0$  otherwise.
We set for $\vec x\in \QQ^{\ify}$, $\beta_0(\vec x)=0$ and
\beq
\beta_k(\vec x):=x_k+\sum_{k<j<\kp}\lan h_{i_k},\al_{i_j}\ran x_j+x_{\kp}.
\label{betak}
\eeq
We define a piecewise-linear operator $S_k=S_{k,\io}$ on $(\QQ^{\ify})^*$ by
\beq
S_k(\vp):=\cases
{\vp-\vp_k\beta_k & if $\vp_k>0$,\cr
 \vp-\vp_k\beta_{\km} & if $\vp_k\leq 0$.\cr}
\label{Sk}
\eeq
Here we set
\beqnn
\Xi_{\io} &:=  &\{S_{j_l}\cd S_{j_2}S_{j_1}x_{j_0}\,|\,
l\geq0,j_0,j_1,\cd,j_l\geq1\},\\
\Sigma_{\io} & := &
\{\vec x\in \ZZ^{\ify}\subset \QQ^{\ify}\,|\,\vp(\vec x)\geq0\,\,{\rm for}\,\,
{\rm any}\,\,\vp\in \Xi_{\io}\}.
\eeqnn
We impose on $\io$ the following positivity assumption:
$$
{\hbox{if $\km=0$ then $\vp_k\geq0$ for any 
$\vp(\vec x)=\sum_k\vp_kx_k\in \Xi_{\io}$}}.
$$
\newtheorem{thm2}[df2]{Theorem}
\begin{thm2}[\cite{NZ}]
Let $\io$ be a sequence of indices satisfying $(\ref{seq-con})$ 
and the positivity
assumption, and  $\Psi_{\io}:B(\ify)\hookrightarrow \ZZ^{\ify}_{\io}$ 
be the corresponding
Kashiwara embedding. Then we have 
${\rm Im}(\Psi_{\io})(\cong B(\ify))=\Sigma_{\io}$.
\end{thm2}

{\sl Remark.\,\,}
We shall show the example 
of the sequence $\io$ which does not satisfy
the positivity assumption.
It will be given in the end of Sect.6.

\subsection{Global crystal base}

In this subsection, we recall several facts 
about global crystal bases
(see \cite{K1},\cite{K6}).

We define a $\QQ$-algebra automorphism $-$ of $\uq$ by:
$\ovl q=q^{-1},$  $\ovl{q^h}=q^{-h}$, $\ovl{e_i}=e_i$, $\ovl{f_i}=f_i$.

Let $\uqq$ be the sub-$\QQ[q,q^{-1}]$-algebra of $\uqm$ generated by
$f^{(n)}_i=f^n_i/[n]_i!$. and $V_{\QQ}(\lm):=
\uqq u_{\lm}$ for $\lm\in P_+$.
Let $p_{\ify}:L(\ify)\rightarrow L(\ify)/qL(\ify)$
(resp. $p_{\lm}:L(\lm)\rightarrow L(\lm)/qL(\lm)$)
be the canonical projection.

\begin{pro2}[\cite{K1}]
\label{thm-triple}
The map $p_{\ify}$ $($resp. $p_{\lm}$$)$ gives rise to the $\QQ$-linear isomorphism:
\beq
\uqq\cap L(\ify)\cap \ovl L(\ify)\mapright{\sim}L(\ify)/qL(\ify)\q
({\rm resp. \,\,}
V_{\QQ}(\lm)\cap L(\lm)\cap \ovl L(\lm)\mapright{\sim}L(\lm)/qL(\lm)).
\label{triple}
\eeq
\end{pro2}
Let us denote  the inverse of this isomorphism by $G$.
The set of inverse image of crystal base 
$\{G(b)\,|\,b\in B(\ify)\,\,({\rm resp.\,\,}B(\lm))\}$
is called {\it global $($crystal$)$ base} of $\uqm$ (resp. $V(\lm)$).
The global base holds the following remarkable property
(\cite[Theorem 7]{K1},\cite[(6.3)]{K6}):
\beq
f_i^n\uqm=\bigoplus_{b\in B(\ify),\,\,\vep_i(b)\geq n}\QQ(q)G(b).
\label{fi-uqm}
\eeq
As we have seen that the anti-automorphism
$*$ preserves $\uqm$, and furthermore, has the
property (\ref{*-sta}), 
which implies that the action of $*$ commutes with $p_{\ify}$ and then
we have $G(b^*)=G(b)^*$.
Thus, applying $*$ on (\ref{fi-uqm}) we obtain,

\beq
\uqm f_i^n=\bigoplus_{b\in B(\ify),\,\,\vep_i(b^*)\geq n}\QQ(q)G(b).
\label{uqm-fi}
\eeq
For a dominant integral weight $\lm$, let $\pi_{\lm}$ be the projection
$\uqm\rightarrow V(\lm)$ as in 2.2. By (\ref{vlm}) we know that
\beq
{\rm Ker}(\pi_{\lm})=\sum_{i}\uqm f_i^{1+\lan h_i,\lm\ran}.
\label{kernel}
\eeq
By virtue of (\ref{uqm-fi}) and (\ref{kernel}) we have:
\beq
\pi_{\lm}(G(b))=0\q \Longleftrightarrow\q  \vep_i(b^*)>\lan h_i,\lm\ran
 \,\,{\rm for \,\,some\,\,}i\in I.
\label{pi=0}
\eeq

\begin{pro2}
\label{e*}
For $b\in B(\ify)$ and $\lm\in P_+$,
$\what\pi_{\lm}(b)\ne0$ if and only if $\vep^*_i(b)\leq \lan h_i,\lm\ran$
for any $i\in I$, where we set $\vep^*_i(b):=\vep_i(b^*)$.
\end{pro2}

\vskip5pt
{\sl Proof.\,\,}
Since $p_{\lm}\circ \pi_{\lm}=\wpi\circ p_{\ify}$
(see \cite{K1}), it follows from 
Proposition \ref{thm-triple} and (\ref{pi=0})
that
$\pi_{\lm}(G(b))=0\Longleftrightarrow \what\pi_{\lm}(b)=0$.
Thus, we get the desired result.\qed


\section{Embedding of $B(\lm)$}
\setcounter{equation}{0}
\renewcommand{\thesection}{\arabic{section}}
\renewcommand{\theequation}{\thesection.\arabic{equation}}

In this section, $\lm$ is supposed to be a dominant integral weight.
\subsection{Morphisms of Crystals}

We shall introduce a new morphism of crystals by modifying 
the map $\wpi$.
Let $R_{\lm}$ be the crystal defined in Example \ref{Example:crystal} (ii).
Consider the crystal $B(\ify)\ot R_{\lm}$ and  define the map
\beq
\Phi_{\lm}:(B(\ify)\ot R_{\lm})\sqcup\{0\}\longrightarrow B(\lm)\sqcup\{0\},
\label{philm}
\eeq
by $\Phi_{\lm}(0)=0$ and $\Phi_{\lm}(b\ot r_{\lm})=\wpi(b)$ for $b\in B(\ify)$.
We set
$$
\wtil B(\lm):=
\{b\ot r_{\lm}\in B(\ify)\ot R_{\lm}\,|\,\Phi_{\lm}(b\ot r_{\lm})\ne 0\}.
$$

\newtheorem{thm3}{Theorem}[section]
\begin{thm3}
\label{ify-lm}
\begin{enumerate}
\item
The map $\Phi_{\lm}$ becomes a surjective strict morphism of crystals
$B(\ify)\ot R_{\lm}\longrightarrow B(\lm)$.
\item
$\wtil B(\lm)$ is a subcrystal of $B(\ify)\ot R_{\lm}$, 
and $\Phi_{\lm}$ induces the
isomorphism of crystals $\wtil B(\lm)\mapright{\sim} B(\lm)$.
\item
We have
\beq
\wtil B(\lm)=\{b\ot r_{\lm}\in B(\ify)\ot R_{\lm}\,|\,
\vep^*_i(b)\leq \lan h_i,\lm\ran\,\,{\rm for }\,\,{\rm any}\,\,i\in I\}.
\label{tilblm}
\eeq
\end{enumerate}
\end{thm3}
The proof of this theorem will be given in the next subsection.

Let us denote $\ZZ^{\ify}_{\io}\ot R_{\lm}$ by 
$\ZZ^{\ify}_{\io}[\lm]$. Here note that
since the crystal $R_{\lm}$ has  only one element,
as a set we can identify $\ZZ^{\ify}_{\io}[\lm]$ with
$\ZZ^{\ify}_{\io}$ but their crystal structures are different.
By Theorem \ref{ify-lm}, we have the strict embedding of crystals
(see also \cite{J}):
$$
\Omega_{\lm}:B(\lm)(\cong \wtil B(\lm))\hookrightarrow B(\ify)\ot R_{\lm}.
$$
Combining $\Omega_{\lm}$ and the
Kashiwara embedding $\Psi_{\io}$,
we obtain the following:
\begin{thm3}
\label{embedding}
There exists the unique  strict embedding of crystals
\beq
\Psi_{\io}^{(\lm)}:B(\lm)\stackrel{\Omega_{\lm}}{\hookrightarrow}
B(\ify)\ot R_{\lm}
\stackrel{\Psi_{\io}\ot {\rm id}}{\hookrightarrow}
\ZZ^{\ify}_{\io}\ot R_{\lm}=:\ZZ^{\ify}_{\io}[\lm],
\label{Psi-lm}
\eeq
such that $\Psi^{(\lm)}_{\io}(u_{\lm})=(\cd,0,0,0)\ot r_{\lm}$.
\end{thm3}

\vskip5pt

The main result of the present paper is an explicit 
description of the image of $\Psi^{(\lm)}_{\io}$
$(\cong \wtil B(\lm))$ as a part in
$B(\ify)\ot R_{\lm}\hookrightarrow \ZZ^{\ify}_{\io}[\lm]$,
which will be given in Sect.4.

\subsection{Proof of Theorem \ref{ify-lm} }

\nd
Before showing Theorem \ref{ify-lm},
we see the following lemmas:
\newtheorem{lm3}[thm3]{Lemma}
\begin{lm3}
\label{ep*}
For $b\in B(\ify)$, suppose that $\til e_i b\ne0$. Then we have
$$
\vep^*_j(\eit  b)= \vep^*_j(b) \,\,(i\ne j)\,\,\,
{\rm and}\,\,\,
\vep^*_i(\til e_i b)\leq \vep^*_i(b).
$$
\end{lm3}

{\sl Proof.\,\,}
Let $\Psi_j:B(\ify)\hookrightarrow B(\ify)\ot B_j$ 
$(u_{\ify}\mapsto u_{\ify}\ot (0)_j)$ be the 
strict embedding as in (\ref{Psii}),
which satisfies that for $b\in B(\ify)$
$\Psi_i(b)=b_1\ot \til f_j^m(0)_j$ where $m=\vep^*_j(b)$ and
$b_1=(\til e_j^m b^*)^*$.
If $i\ne j$, $\Psi_j(\til e_i b)=\til e_i\Psi_j(b)=
(\til e_i b_1)\ot \til f_j^m(0)_j$ by (\ref{tensor-e}) and then
we have $\vep^*_j(\eit b)=m=\vep^*_j(b)$.
In the case $i=j$, we have
$\eit(b_1\ot \fit^m(0)_i)=\eit b_1\ot \fit^m(0)_i$ or
$b_1\ot \fit^{m-1}(0)_i$ ($m\geq 1$). This implies that
$\vep^*_i(\eit b)\leq \vep^*_i(b)$.\qed

\begin{lm3}
\label{eitu=eitb}
Suppose that $\wpi(b)\ne0$ for $b\in B(\ify)$. Then 
$\eit \wpi(b)=0$ if and only if $\eit b=0$.
\end{lm3}

{\sl Proof.}\,\,
We assume $\eit \wpi(b)=0$.
Since $\eit \wpi(b)=\wpi(\eit b)$ if $\wpi(b)\ne0$ by (\ref{epi=pie}),
 we have $\wpi(\eit b)=0$.
If $\eit b\ne0$, it follows from Lemma \ref{ep*} that  for any $j\in I$
$\vep^*_j(\eit b)\leq \vep^*_j(b)\leq \lan h_j,\lm\ran$,
which contradicts $\wpi(\eit b)=0$ by Proposition \ref{e*}.
Hence, we have $\eit b=0$.
On the other hand, it is trivial that if $\eit b=0$, then
$\eit\wpi(b)=\wpi(\eit b)=0$ by (\ref{epi=pie}).\qed

\vskip5pt

{\sl Proof of Theorem \ref{ify-lm}.}
The statement (iii) of the theorem is an 
immediate consequence of Proposition \ref{e*}.

Let us show (i).
The surjectivity follows from the one for the map $\wpi$.
So we try to prove that $\Phi_{\lm}$ is a strict morphism of crystals. 
To do this, according to Definition \ref{df:mor} (i)
it suffices to show: for $u\in B(\ify)\ot R_{\lm}$,

\renewcommand{\labelenumi}{(\arabic{enumi})}
\begin{enumerate}
\item
$wt(\Phi_{\lm}(u))=wt(u)$ if $\Phi_{\lm}(u)\ne0$,
\item
$\vep_i(\Phi_{\lm}(u))=\vep_i(u)$ for any $i$ if $\Phi_{\lm}(u)\ne0$,
\item
$\vp_i(\Phi_{\lm}(u))=\vp_i(u)$ for any $i$ if $\Phi_{\lm}(u)\ne0$,
\item
$\eit \Phi_{\lm}(u)=\Phi_{\lm}(\eit u)$ for any $i$,
\item
$\fit \Phi_{\lm}(u)=\Phi_{\lm}(\fit u)$ for any $i$.
\end{enumerate}

\renewcommand{\labelenumi}{(\roman{enumi})}

Let us show (1).
For $u=b\ot \rlm=
(\til f_{i_l}\cd\til f_{i_1}u_{\ify})\ot r_{\lm}\in B(\ify)\ot R_{\lm}$,
we have $\Phi_{\lm}(u)=\wpi(b)=\til f_{i_l}\cd\til f_{i_1}u_{\lm}$
since any $\fit$ commutes with $\wpi$.
It follows immediately that
$wt(u)=\lm-\al_{i_1}-\cd-\al_{i_l}
=wt(\Phi_{\lm}(u))$ if $\Phi_{\lm}(u)\ne0$.

\vskip5pt

In the case $\wpi(b)\ne0$,
it follows from Lemma \ref{eitu=eitb} that
$\eit b=0$ if and only if $\eit \wpi(b)=0$.
This means if $\wpi(b)\ne0$,
\beq
\vep_i(b)=\vep_i(\wpi(b)).
\label{vep=vep}
\eeq
Furthermore, by (1), (\ref{vp=vep+wt}), 
(\ref{vep-vp}) and (\ref{vep=vep}) we have
\beq
0\leq \vp_i(\wpi(b))
=\lan h_i,\lm\ran+\vp_i(b).
\label{vp(u)}
\eeq
It follows from (\ref{tensor-vep}), (\ref{vep=vep}) and (\ref{vp(u)})
that 
\beqnn
\vep_i(b\ot r_{\lm})
&=&{\rm max}(\vep_i(b),\vep_i(r_{\lm})-\lan h_i,wt(b)\ran)
={\rm max}(\vep_i(b),\vep_i(b)-\vp_i(b)-\lan h_i,\lm\ran)\\
&=&\vep_i(b)=\vep_i(\wpi(b))=\vep_i(\Phi_{\lm}(u)).
\eeqnn
Now we obtained (2).

The statement (3) is derived immediately from (1), (2) and
(\ref{vp=vep+wt}).

%

Let us show (4), namely,
$\eit\Phi_{\lm}(b\ot\rlm)=\Phi_{\lm}(\eit(b\ot\rlm))$ for any $i$.

\nd
First, for $u=b\ot r_{\lm}$ suppose $\Phi_{\lm}(u)=\wpi(b)\ne0$.
Then by (\ref{epi=pie}), we have $\eit\Phi_{\lm}(u)=\eit\wpi(b)=\wpi(\eit b)$.
Thus, it suffices to show
\beq
\eit(b\ot r_{\lm})=(\eit b)\ot r_{\lm}.
\label{eiteit}
\eeq
By (\ref{vp(u)}) we have
$\vp_i(b)\geq -\lan h_i,\lm\ran=\vep_i(r_{\lm})$, which means
(\ref{eiteit}) by (\ref{tensor-e}).

Next, we consider the case
$\Phi_{\lm}(u)=\wpi(b)=0.$
It sufficies to show
\beq
\Phi_{\lm}(\eit(b\ot r_{\lm}))=0.
\label{eit-0}
\eeq
If $\eit(b\ot r_{\lm})=0$, there is nothing to show.
So we consider the case $\eit(b\ot r_{\lm})\ne0$.
Since $\eit r_{\lm}=0$, we have $\eit(b\ot r_{\lm})=(\eit b)\ot r_{\lm},$
which implies 
\beq
\vp_i(b)\geq \vep_i(r_{\lm})=-\lan h_i,\lm\ran,
\label{vp1}
\eeq
by (\ref{tensor-e}).
Now, assuming 
$\Phi_{\lm}((\eit b)\ot r_{\lm})=\wpi(\eit b)\ne0,$
we shall derive a contradiction.
We have $\fit(\wpi(\eit b))=\wpi(\fit\eit b)=\wpi(b)=0$.
Since $\wpi(\eit b)\ne0$, we obtain
$\vp_i(\wpi(\eit b))=0$ (see (\ref{vep-vp})).
Thus, taking into account (\ref{vp=vep+wt}), (\ref{+alpha}), 
(\ref{vep-vp}) (\ref{vep-vp-ify}) and (\ref{vep=vep}), we have
\beqnn
0 &=& \vp_i(\wpi(\eit b))=
\lan h_i,wt(\wpi(\eit b))\ran +\vep_i(\wpi(\eit b))
=\lan h_i,\lm+wt(b)+\al_i\ran+\vep_i(\eit b)\\
&=& \lan h_i,\lm\ran +\lan h_i,wt(b)\ran +2+\vep_i(b)-1
=\lan  h_i,\lm\ran +\vp_i(b)-\vep_i(b)+\vep_i(b)+1\\
&=&\lan h_i,\lm\ran+\vp_i(b)+1.
\eeqnn
Thus, we have 
$\vp_i(b)=-\lan h_i,\lm\ran-1=\vep_i(r_{\lm})-1<\vep_i(r_{\lm})$,
which contradicts (\ref{vp1}). 
Therefore, we have $\eit(b\ot r_{\lm})=0$ and completed to prove (4).

Finally, let us show (5).
Since $\fit$ commutes with $\wpi$, 
for $u=b\ot r_{\lm}$ we have
$\fit\Phi_{\lm}(u)=\fit\wpi(b)=\wpi(\fit b).$
Thus, if $\fit u=(\fit b)\ot r_{\lm}$, 
we have
$\wpi(\fit b)=\Phi_{\lm}(\fit u)$, which means (5).
So we consider the case $\fit u=b\ot \fit r_{\lm}=0$.
In this case, we shall try to show that $\fit\Phi_{\lm}(u)=0$.
It follows from (\ref{tensor-f})
that we have
\beq
\vp_i(b)\leq \vep_i(r_{\lm})=-\lan h_i,\lm\ran.
\label{vp-con}
\eeq
If $\Phi_{\lm}(u)=\wpi(b)=0$, 
there is nothing to show. So we may consider the case
$\Phi_{\lm}(u)=\wpi(b)\ne 0$.
Assuming 
\beq
\fit\Phi_{\lm}(u)=\fit\wpi(b)\ne0, 
\label{fpibne0}
\eeq
we shall derive a contradiction. 
The assumption (\ref{fpibne0}) means $\vp_i(\wpi(b))>0$.
Thus, by (\ref{vp(u)})
we have $0< \lan h_i,\lm\ran+\vp_i(b)$. It 
implies $\vp_i(b)>-\lan h_i,\lm\ran=\vep_i(r_{\lm})$,
which contradicts (\ref{vp-con}).
Now we obtain $\fit\Phi_{\lm}(u)=0$.
Thus, we have completed to show (5) and then (i).

Let us show (ii).
The condition $\Phi_{\lm}(b\ot r_{\lm})\ne0$ is equivalent to
$\wpi(b)\ne0$.
Thus, the map $\phi_{\lm}:=
{\Phi_{\lm}|}_{\wtil B(\lm)}:\wtil B(\lm)\longrightarrow 
B(\lm)$ is bijective by (\ref{wpi-bij}).
So if we show that $\wtil B(\lm)$ is 
stable by the actions of $\eit$ and $\fit$,
it follows from (i) that the map $\phi_{\lm}$ is a 
strict morphism of crystals.
Let us see the stability of $\wtil B(\lm)$, namely, 
that if $\Phi_{\lm}(\eit(b\ot r_{\lm}))=0$
(resp. $\Phi_{\lm}(\fit(b\ot r_{\lm}))=0$)
for $b\ot r_{\lm}\in \wtil B(\lm)$,
then $\eit(b\ot r_{\lm})=0$ (resp. $\fit(b\ot r_{\lm})=0$).

First, for $b\ot r_{\lm}\in \wtil B(\lm)$ 
suppose that $\Phi_{\lm}(\eit(b\ot r_{\lm}))=0$,
which implies $\eit\wpi(b)=0$ by (i). 
Since $\wpi(b)\ne0$, we have $\eit b=0$ by Lemma \ref{eitu=eitb}.
Thus, we obtain $\eit(b\ot r_{\lm})=0$ in view of (\ref{tensor-e})
and $\eit(r_{\lm})=0$.
Next, for $b\ot r_{\lm}\in \wtil B(\lm)$ 
suppose that $\Phi_{\lm}(\fit(b\ot r_{\lm}))=0$,
which implies $\fit\wpi(b)=0$.
It follows from $\wpi(b)\ne0$ that 
$\vp_i(\wpi(b))=0$. Then, 
we have
$$
\vp_i(\wpi(b))=\vp_i(b)+\lan h_i,\lm\ran=0\,\,
\Longleftrightarrow \vp_i(b)=-\lan h_i,\lm\ran =\vep_i(r_{\lm}).
$$
This means 
$\fit(b\ot r_{\lm})=b\ot \fit r_{\lm}=0$.
Now we have completed to prove (ii) and then Theorem \ref{ify-lm}.
\qed

\newtheorem{ex3}[thm3]{Example}
\begin{ex3}
Let us see the simplest example $\ge=\ssl_2$-case.
Let $u_{\ify}$ be the highest weight vector in $B(\ify)$. 
Then we have 
$B(\ify)=\{\ft^n u_{\ify}\}$.
The crystal graph of $B(\ify)$ is as follows:

\nd
\unitlength 1pt
\begin{picture}(800,15)
\put(0,0){\circle{10}}
\put(0,0){\makebox(0,0){0}}
\put(5,0){\vector(1,0){30}}
\put(40,0){\circle{10}}
\put(40,0){\makebox(0,0){1}}
\put(45,0){\vector(1,0){30}}
\put(80,0){\circle{10}}
\put(80,0){\makebox(0,0){2}}
\put(85,0){\vector(1,0){30}}
\put(120,0){\circle{10}}
\put(125,0){\vector(1,0){30}}
\put(160,0){\circle{10}}
\put(165,0){\vector(1,0){30}}
\put(200,0){\circle{10}}
\put(200,0){\makebox(0,0){m}}
\put(205,0){\vector(1,0){30}}
\put(240,0){\circle{10}}
\put(245,0){\vector(1,0){30}}
\put(280,0){\circle{10}}
\put(285,0){\vector(1,0){30}}
\put(322,0){\line(1,0){3}}
\put(328,0){\line(1,0){3}}
\put(334,0){\line(1,0){3}}
\put(340,0){\line(1,0){3}}
\end{picture}
\vskip10pt

\nd
where 
\unitlength 1pt
\begin{picture}(10,5)
\put(4,4){\circle{10}}
\put(4,4){\makebox(0,0){x}}
\end{picture}$=\ft^x u_{\ify}$.

Next, let us see the crystal grah of $B(\ify)\ot R_m$ 
$(m\in\ZZ_{\geq0})$.
We know that $\vp(\ft^n u_{\ify})=-n$ and $\vep(r_m)=-m$. 
Then, by (\ref{tensor-f}) we have 
$$
\ft(\ft^n u_{\ify}\ot r_m)=
\cases{\ft^{n+1}u_{\ify}\ot r_m& if $n<m$, \cr
       \ft^n u_{\ify}\ot\ft(r_m)=0 & if $n \geq m$.\cr}
$$
Thus, the crystal graph of $B(\ify)\ot R_m$ is:

\nd
\unitlength 1pt
\begin{picture}(800,15)
\put(-5,-5){\framebox(10,10){0}}
\put(5,0){\vector(1,0){30}}
\put(35,-5){\framebox(10,10){1}}
\put(45,0){\vector(1,0){30}}
\put(75,-5){\framebox(10,10){2}}
\put(85,0){\vector(1,0){30}}
\put(115,-5){\framebox(10,10){}}
\put(125,0){\vector(1,0){30}}
\put(155,-5){\framebox(10,10){}}
\put(165,0){\vector(1,0){30}}
\put(195,-5){\framebox(10,10){m}}
\put(235,-5){\framebox(10,10){}}
\put(275,-5){\framebox(10,10){}}
\put(322,0){\line(1,0){3}}
\put(328,0){\line(1,0){3}}
\put(334,0){\line(1,0){3}}
\put(340,0){\line(1,0){3}}
\end{picture}
\vskip10pt

\nd
where \framebox{x}$=\ft^x u_{\ify}\ot r_m$.
The connected component including 
$\framebox{0}=u_{\ify}\ot r_m$ is isomorphic
to the crystal $B(m)$ associated with the $m+1$-dimensional 
irreducible module $V(m)$.
In the subsequent section, we shall see how to remove the 
vectors cut off from $B(\lm)$ $($in this case, the vectors 
$\{\framebox{x}|x> m\}$.$)$.
\end{ex3}

\section{Polyhedral Realization of $B(\lm)$}
\setcounter{equation}{0}
\renewcommand{\thesection}{\arabic{section}}
\renewcommand{\theequation}{\thesection.\arabic{equation}}

\subsection{Crystal structure of $\ZZ^{\ify}[\lm]$}

We shall give an explicit crystal structure of
$\ZZ^{\ify}[\lm]$ in a similar manner to \cite{NZ}.
Fix a sequence of indices $\io:=(i_k)_{k\geq 1}$ satisfying the condition
(\ref{seq-con}) and a weight $\lm\in P$.
(In this subsection, we do not necessarily assume that 
$\lm$ is dominant.)
As we stated in 3.1, 
we can identify $\ZZ^{\ify}$ with $\ZZ^{\ify}[\lm]$
as a set. Thus $\ZZ^{\ify}[\lm]$ can be regarded 
as a subset of $\QQ^{\ify}$, and then 
we denote an element in $\ZZ^{\ify}[\lm]$
by $\vec x=(\cd,x_k,\cd,x_2,x_1)$.
For $\vec x=(\cd,x_k,\cd,x_2,x_1)\in \QQ^{\ify}$
we define the linear functions
\beqn
\sigma_k(\vec x)&:= &x_k+\sum_{j>k}\lan h_{i_k},\al_{i_j}\ran x_j,
\q(k\geq1)
\label{sigma}\\
\sigma_0^{(i)}(\vec x)
&:= &-\lan h_i,\lm\ran+\sum_{j\geq1}\lan h_i,\al_{i_j}\ran x_j,
\q(i\in I)
\label{sigma0}
\eeqn
Here note that
since $x_j=0$ for $j\gg0$ on $\QQ^{\ify}$,
the functions $\sigma_k$ and $\sigma^{(i)}_0$ are
well-defined.
Let $\sigma^{(i)} (\vec x)
 := {\rm max}_{k: i_k = i}\sigma_k (\vec x)$, and
\beq
M^{(i)} = M^{(i)} (\vec x) :=
\{k: i_k = i, \sigma_k (\vec x) = \sigma^{(i)}(\vec x)\}.
\label{m(i)}
\eeq
Note that
$\sigma^{(i)} (\vec x)\geq 0$, and that
$M^{(i)} = M^{(i)} (\vec x)$ is a finite set
if and only if $\sigma^{(i)} (\vec x) > 0$.
Now we define the maps
$\eit: \ZZ^{\ify}[\lm] \sqcup\{0\}\lar \ZZ^{\ify}[\lm] \sqcup\{0\}$
and
$\fit: \ZZ^{\ify}[\lm] \sqcup\{0\}\lar \ZZ^{\ify}[\lm] \sqcup\{0\}$ 
by setting $\eit(0)=\fit(0)=0$, and 
\beq
(\fit(\vec x))_k  = x_k + \delta_{k,{\rm min}\,M^{(i)}}
\,\,{\rm if }\,\,\sigma^{(i)}(\vec x)>\sigma^{(i)}_0(\vec x);
\,\,{\rm otherwise}\,\,\fit(\vec x)=0,
\label{action-f}
\eeq
\beq
(\eit(\vec x))_k  = x_k - \delta_{k,{\rm max}\,M^{(i)}} \,\, {\rm if}\,\,
\sigma^{(i)} (\vec x) > 0\,\,
{\rm and}\,\,\sigma^{(i)}(\vec x)\geq\sigma^{(i)}_0(\vec x) ; \,\,
 {\rm otherwise} \,\, \eit(\vec x)=0,
\label{action-e}
\eeq
where $\del_{i,j}$ is the Kronecker's delta.
We also define the weight function and the functions
$\vep_i$ and $\vp_i$ on $\ZZ^{\ify}[\lm]$ by
\beq
\begin{array}{l}
wt(\vec x) :=\lm -\sum_{j=1}^{\ify} x_j \al_{i_j}, \,\,
\vep_i (\vec x) := {\rm max}(\sigma^{(i)} (\vec x),\sigma^{(i)}_0(\vec x))\\
\vp_i (\vec x) := \lan h_i, wt(\vec x) \ran + \vep_i(\vec x).
\end{array}
\label{wt-vep-vp}
\eeq
We will denote this crystal by $\ZZ^{\ify}_{\io}[\lm]$.
Note that, in general, the subset $\ZZ^{\ify}_{\geq 0}[\lm]$ is not
a subcrystal of $\ZZ^{\ify}_{\io}[\lm]$ since it is not
stable under the action of $\eit$'s.

\subsection{The image of $\Psi^{(\lm)}_{\io}$}


As in the previous sections, we fix a sequence of indices $\io$
satisfying (\ref{seq-con}) and take a dominant integral weight 
$\lm\in P_+$.
For $k\geq1$ let $k^{(\pm)}$ be  the ones in 2.4.
Let $\beta_k^{(\pm)}(\vec x)$ be linear functions given by
\beqn
\beta_k^{(+)} (\vec x) & = & \sigma_k (\vec x) - \sigma_{\kp} (\vec x),
\label{beta}\\
\beta_k^{(-)} (\vec x) & = &
\cases{
\sigma_{\km} (\vec x) - \sigma_k (\vec x)& if $\km>0$,\cr
\sigma_0^{(i_k)} (\vec x) - \sigma_k (\vec x)& if $\km=0$,\cr}
\label{beta--}
\eeqn
where the functions $\sigma_k$ and $\sigma^{(i)}_0$
are defined by (\ref{sigma}) and (\ref{sigma0}).
Since $\lan h_{i},\al_{i}\ran   = 2$ for any $i \in I$, we have
\beqn
\beta_k^{(+)}(\vec{ x}) & = &  x_k+\sum_{k<j<\kp}
\lan h_{i_k},\al_{i_j}\ran x_j+x_{\kp},
\label{beta+}\\
\beta_k^{(-)}(\vec{ x}) & = &
\cases{
x_{\km}+\sum_{\km<j<k}\lan h_{i_k},\al_{i_j}\ran x_j+x_k & if $\km>0$,\cr
-\lan h_{i_k},\lm\ran
+\sum_{1\leq j<k}\lan h_{i_k},\al_{i_j}\ran x_j+x_k & if $\km=0$.\cr}
\label{beta-}
\eeqn
Here note that
$$
\beta_k^{(+)}=\beta_k,\qq
\beta_k^{(-)}=\beta_{\km}  {\hbox{ \,\,if\,\, $\km>0$}}.
$$

Using this notation, for every $k \geq 1$, we define 
an operator
$\what S_k = \what S_{k,\io}$ for a linear function $\vp(\vec x)=c+\sum_{k\geq 1}\vp_kx_k$
$(c,\vp_k\in\QQ)$ on $\QQ^{\ify}$ by:

\beq
\what S_k\,(\vp) :=\cases{
\vp - \vp_k \beta_k^{(+)} & if $\vp_k > 0$,\cr
\vp - \vp_k \beta_k^{(-)} & if $\vp_k \leq 0$.\cr}
\label{S_k}
\eeq
An easy check shows that
$(\what S_k)^2=\what S_k.$

For the fixed sequence $\io=(i_k)$, 
in case $\km=0$ for $k\geq1$, there exists unique $i\in I$ such that $i_k=i$.
We denote such $k$ by $\io^{(i)}$, namely, $\io^{(i)}$ is the first
number $k$ such that $i_k=i$.

Here we set
\beq
\lm^{(i)}(\vec x):=
-\beta^{(-)}_{\io^{(i)}}(\vec x)=\lan h_i,\lm\ran-\sum_{1\leq j<\io^{(i)}}
\lan h_i,\al_{i_j}\ran x_j-x_{\io^{(i)}}.
\label{lmi}
\eeq

For $\io$ and a dominant integral weight $\lm$,
let $\Xi_{\io}[\lm]$ be the set of all linear functions
generatd by applying $\what S_k=\what S_{k,\io}$ 
on the functions $x_j$ ($j\geq1$)
and $\lm^{(i)}$ ($i\in I$), namely,
\beq
\begin{array}{ll}
\Xi_{\io}[\lm]&:=\{\what S_{j_l}\cd\what S_{j_1}x_{j_0}\,
:\,l\geq0,\,j_0,\cd,j_l\geq1\}
\\
&\cup\{\what S_{j_k}\cd \what S_{j_1}\lm^{(i)}(\vec x)\,
:\,k\geq0,\,i\in I,\,j_1,\cd,j_k\geq1\}.
\end{array}
\label{Xi}
\eeq
Now we set
\begin{equation}
\Sigma_{\io}[\lm]
:=\{\vec x\in \ZZ^{\ify}_{\io}[\lm](\subset \QQ^{\ify})\,:\,
\vp(\vec x)\geq 0\,\,{\rm for \,\,any }\,\,\vp\in \Xi_{\io}[\lm]\}.
\label{Sigma}
\end{equation}

For a sequence $\io$ and a domiant integral weight $\lm$, a pair
$(\io,\lm)$ is called {\it ample}
if $\Sigma_{\io}[\lm]\ni\vec 0=(\cd,0,0)$.

\newtheorem{thm4}{Theorem}[section]
\begin{thm4}
\label{main}
Suppose that $(\io,\lm)$ is ample.
Let $\Psi^{(\lm)}_{\io}:B(\lm)\hookrightarrow \ZZ^{\ify}_{\io}[\lm]$
be the embedding as in (\ref{Psi-lm}). Then the image
${\rm Im}(\plm)(\cong B(\lm))$ is equal to $\Sigma_{\io}[\lm]$.
\end{thm4}

{\sl Proof.\,\,}
Taking into account of (\ref{def-blm}) and Theorem \ref{embedding},
the image ${\rm Im}(\plm)$ is a subcrystal of $\ZZ^{\ify}[\lm]$
obtained by
applying $\fit$'s to $\plm(u_{\lm})=\vec 0=(\cd,0,0)$, that is,
\beq
{\rm Im}(\plm)=
\{\til f_{i_l}\cdot\cdot \til f_{i_1}
\plm(u_{\lm})\,|\,i_j\in I,l\geq 0\}\setminus\{0\}.
\eeq
By the explicit description of $\fit$ in (\ref{action-f}),
we know that ${\rm Im}(\plm)\subset \ZZ^{\ify}_{\geq 0}[\lm]$.
Since the pair $(\io,\lm)$ is ample, $\Sigma_{\io}[\lm]\ni\vec 0$.
Thus, the inclusion ${\rm Im}(\plm)\subset \Sigma_{\io}[\lm]$ follows
from the fact  that the set $\Sigma_{\io}[\lm]$ is closed by the
actions of $\fit$'s, namely,
$\fit \Sigma_{\io}[\lm]\subset \Sigma_{\io}[\lm]\sqcup\{0\}$ for any $i\in I$.
Let us show this.
For $\vec x=(\cd, x_2,x_1)\in\Sigma_{\io}[\lm]$ and
$i\in I$, suppose that
$\fit \vec x=(\cd,x_k+1,\cd,x_2,x_1)$ (note that $i_k=i$).
We shall show
\beq
\vp(\fit \vec x)\geq0,
\label{vp(f)}
\eeq
for any $\vp(\vec x)=c+\sum\vp_jx_j\in \Xi_{\io}[\lm]$.
Since $\vp(\fit \vec x)=\vp(\vec x)+\vp_k\geq \vp_k$,
it suffices to consider the case $\vp_k<0$.
By the definition of $M^{(i)}(\vec x)$ in (\ref{m(i)}), 
we know that $k$ is the minimum in $M^{(i)}(\vec x)$.
Thus, it follows from (\ref{action-f}) that 
$\sigma_k(\vec x)>\sigma_{\km}(\vec x)$ if $\km>0$ or 
$\sigma_k(\vec x)>\sigma^{(i)}_0(\vec x)$ if $\km=0$.
Thus, by (\ref{beta--}) we have
$\beta^{(-)}_k(\vec x)<0$. 
Therefore, since the function $\beta^{(-)}_k$ takes an integer value for 
$\vec x\in \ZZ^{\ify}$, 
\beq
\beta^{(-)}_k(\vec x)\leq -1.
\label{beta-<-1}
\eeq
It follows from $\what S_k\vp\in \Xi_{\io}[\lm]$ and 
$\vp_k<0$ that 
$$
\vp(\fit \vec x) = \vp(\vec x)+\vp_k
\geq \vp(\vec x)-\vp_k\beta^{(-)}_k(\vec x)
= (\what S_k\vp)(\vec x)\geq0.
$$
Therefore, we get the inclusion ${\rm Im}(\plm)\subset \Sigma_{\io}[\lm]$.

Let us show the reverse inclusion
$\Sigma_{\io}[\lm]\subset {\rm Im}(\plm)$.
We first show that $\Sigma_{\io}[\lm]$ is
a subcrystal of $\ZZ^{\ify}_{\io}[\lm]$.
Since we have already shown that
$\fit \Sigma_{\io}[\lm]\subset \Sigma_{\io}[\lm]\sqcup\{0\}$ for any $i\in I$,
it is enough to prove that
$\eit \Sigma_{\io}[\lm]\subset \Sigma_{\io}[\lm]\sqcup\{0\}$ for any $i\in I$.
For $\vec x=(\cd,x_2,x_1)\in \Sigma_{\io}[\lm]$ and $i\in I$,
suppose that $\eit\vec x=(\cd,x_k-1,\cd,x_2,x_1)$,
here note that $i_k=i$.
We have to show
\beq
\vp(\eit \vec x)\geq 0.
\eeq
Since $\vp(\eit \vec x)=\vp(\vec x)-\vp_k\geq -\vp_k$,
it suffices to consider the case $\vp_k>0$.
Arguing similarly to the $\fit$ case, 
by (\ref{m(i)}), (\ref{action-e}) and (\ref{beta}), we have
\beq
\beta^{(+)}_k(\vec x)\geq 1.
\label{beta+>1}
\eeq
It follows from $\what S_k\vp\in \Xi_{\io}[\lm]$ and $\vp_k>0$ that
$$
\vp(\eit \vec x) = \vp(\vec x)-\vp_k
 \geq  \vp(\vec x)-\vp_k\beta^{(+)}_k(\vec x)
= (\what S_k\vp)(\vec x)\geq0.
$$

Since $\Sigma_{\io}[\lm]$ is included in 
$\subset\ZZ^{\ify}_{\geq0}[\lm]$ and
is closed under the actions of $\eit$,
for any $\vec x\in \Sigma_{\io}[\lm]$ there exists
$l\gg0$ such that
$\til e_{i_1}\til e_{i_2}\cd\til e_{i_l}\vec x=0$
for  any $i_1,\cd,i_l\in I$.
(Indeed, we can take $l=\sum_j x_j+1$.).
Therefore, to complete the proof,
it is enough to show that
if $\vec x\in \Sigma_{\io}[\lm]$ satisfies
$\eit \vec x=0$ for any $i\in I$, then $\vec x=\vec 0=(\cd,0,0)$.
Indeed, this implies that for any $\vec x\in \Sigma_{\io}[\lm]$
there exist $i_1,i_2,\cd,i_k\in I$ such that
$\vec 0=\til e_{i_1}\til e_{i_2}\cd\til e_{i_k}\vec x$ or
equivalently,
$\vec x=\til f_{i_k}\cd\til f_{i_2}\til f_{i_1}\vec 0$.
It follows that $\Sigma_{\io}[\lm]\subset {\rm Im}(\plm)$.
Now, suppose that
$\vec x\in \Sigma_{\io}[\lm]$ satisfies $\eit \vec x=0$ for any $i\in I$ and
$\vec x\ne\vec 0$.
By (\ref{action-e}),  for any $i\in I$ we have
\beq
\sigma^{(i)}(\vec x)\leq 0\qq{\rm or}\qq
\sigma^{(i)}(\vec x)<\sigma^{(i)}_0(\vec x).
\label{or}
\eeq
By the assumption $\vec x=(\cd,x_2,x_1)\ne\vec 0$, there exists
$j\geq1$ such that $x_j>0$ and $x_k=0$ for $k>j$.
Thus we have $\sigma_j(\vec x)=x_j>0$. This implies
that $\sigma^{(i_j)}(\vec x)\geq x_j>0$. Then there is no
possibility of the first case in (\ref{or}) for 
$i=i_j$.
Now we suppose that
$\sigma^{(i)}(\vec x)<\sigma^{(i)}_0(\vec x)$ for this $i=i_j$.
Hence, we have 
$0<\sigma^{(i)}_0(\vec x)-\sigma^{(i)}(\vec x)\leq
\sigma^{(i)}_0(\vec x)-\sigma_{\io^{(i)}}(\vec x)
=\beta^{(-)}_{\io^{(i)}}(\vec x)$ and
then $\lm^{(i)}(\vec x)(=-\beta^{(-)}_{\io^{(i)}}(\vec x))<0$ 
(see (\ref{lmi})),
which contradicts the definition of $\Sigma_{\io}[\lm]$
in (\ref{Sigma}).
Now we have completed the proof of Theorem \ref{main}
\qed

\subsection{Formula for $\vep^*_i$ and polyhedral realization of $B(\lm)$}

It is important to get the explicit form of 
$\vep^*_i$ in the sense of Proposition \ref{e*}. 
But a direct computation of $\vep^*_i$ seems to be difficult.
So we apply Theorem \ref{main} to calculate the value $\vep^*_i(b)$.
We define the linear form $\xii$ $(i\in I)$ on $\QQ^{\ify}$ by
\beq
\xii(\vec x):=-\sum_{1\leq j<\io^{(i)}}\lan h_i,\al_{i_j}
\ran x_j-x_{\io^{(i)}}
=-\lan h_i,\lm\ran + \lm^{(i)}(\vec x)
\label{xii}
\eeq
Let us define the set of linear forms $\Xi_{\io}^{(i)}$ by
\beq
\Xi_{\io}^{(i)}:=\{S_{j_l}\cd S_{j_1}\xii\,|\,
l\geq0,j_1,\cd,j_l\geq1\},
\label{Xii}
\eeq
and set $\Xi_{\io}^{(\ify)}:=\Xi_{\io}$ (see 2.4).
Here we introduce the {\it strict positivity assumption} for $\io$ as follows:
\beq
{\hbox{
if $\km=0$ then $\vp_k\geq0$ for any $\vp=\sum_k\vp_kx_k\in
\left(\bigcup_{j\in I\sqcup\{\ify\}}\Xi_{\io}^{(j)}\right)
\setminus\{\xii\,|\,i\in I\}$}}
\label{st-posi-ass}
\eeq

\nd
{\sl Remark.\,\,}
Since the form $\xii$ has a negative coefficient for $x_{\io^{(i)}}$,
we remove $\xii$ from $\Xi_{\io}^{(i)}$in the definition of the
strict positivity assumption.

\nd
Now we have the following theorem:
\begin{thm4}
\label{thm-e*}
Let $\io$ be a sequence of indices satisfying (\ref{seq-con}) and
the strict positivity assumption, and $\lm$ be a dominant integral weight.
Then for $i\in I$ and $\vec x\in \Sigma_{\io}$ we have
\beq
\vep^*_i(\vec x)={\rm max}\{-\vp(\vec x)\,|\,
\vp\in \Xi_{\io}^{(i)}\}
\label{e*-des}
\eeq
\end{thm4}

\vskip5pt
{\sl Proof of Theorem \ref{thm-e*}.}\,\,
First, let us show that 
the ampleness is always satisfied
under the strict positivity assumption. 
To do this, we see the following lemma:
\newtheorem{lm4}[thm4]{Lemma}
\begin{lm4}
\label{vp=lm+vp}
Under the strict positivity assumption for $\io$, we have
\beq
\what S_{j_l}\cd \what S_{j_1}x_{j_0}=S_{j_l}\cd S_{j_1}x_{j_0},
\label{eq1}
\eeq
for any $l\geq0$, $j_0,\cd,j_l\geq1$,  and
\beq
\what S_{j_l}\cd \what S_{j_1}\lm^{(i)}(\vec x)=\lan h_i,\lm\ran+S_{j_l}\cd S_{j_1}\xii(\vec x),
\label{eq2}
\eeq
for any $l\geq0$, $j_1,\cd,j_l\geq1$ and $i\in I$, if the L.H.S. of (\ref{eq2}) is non-zero.
\end{lm4}

{\sl Proof.\,\,}
First we show (\ref{eq1}).
By the definition of $\beta_k^{(\pm)}$ and $\beta_k$ we know that
if $\km>0$, $\what S_k=S_k$. Furthermore, even if $\km=0$, 
under the positivity assumption, 
$\what S_k=S_k$ because in this case their actions
are given by  using only $\beta^{(+)}_k=\beta_k$.

Next we shall see (\ref{eq2}).
Let us show it by the induction on $l$.
If $l=0$, (\ref{eq2}) is just the equation
$\lm^{(i)}=\lan h_i,\lm\ran+\xii$ in (\ref{lmi}).
Now we assume (\ref{eq2}) for $l>0$ 
and write $\vp(\vec x)=c+\sum_k\vp_kx_k\ne0$ for the both sides of (\ref{eq2}).
First, if $\km\ne0$, we have $\what S_k\vp=S_k \vp$ since $\beta_k^{(+)}=\beta_k$ and
$\beta_k^{(-)}=\beta_{\km}$. If $\km=0$ and $\vp\ne \lm^{(i)}$, by the strict positivity assumption,
we have $\vp_k\geq0$ and then
$$
\what S_k \vp=\vp-\vp_k\beta_k^{(+)}=\vp-\vp_k\beta_k=S_k\vp.
$$
Finally, we consider the case $\km=0$ and $\vp=\lm^{(i)}$.
In this case we have $k=\io^{(j)}$ for some $j\in I$.
By the explicit form of $\lm^{(i)}$
we have $\vp_{\io^{(j)}}=-\lan h_i,\al_{\io^{(j)}}\ran \geq0$
for $j\ne i$ in $\vp=\lm^{(i)}$. Thus, if $k=\io^{(j)}$ ($j\ne i$),
$\what S_k\vp=S_k\vp$ by the fact $\beta^{(+)}_k=\beta_k$.
If $k=\io^{(i)}$, the coefficient of $x_{\io^{(i)}}$ in $\lm^{(i)}$ is $-1$.
In this case, 
$S_k\vp=\vp$ and $\what S_k\vp=0$. Then this is not the case of (\ref{eq2}).
\qed

This lemma implies that under the strict positivity assumption,
any linear function in $\Xi_{\io}[\lm]$ has a non-negative coefficient
$0$ or $\lan h_i,\lm\ran$, which means $(\io,\lm)$ is ample. Therefore,
we have 
\beq
\wtil B(\lm)=\{\vec x\in B(\ify)\ot R_{\lm}\subset\ZZ^{\ify}_{\io}[\lm]\,|\,
\begin{array}{l}
\lan h_i,\lm\ran+\vp(\vec x)\geq0\\
{\rm for \,\,any\,\,}i\in I
\,\,{\rm and \,\,}\vp\in \Xi_{\io}^{(i)}\}.
\end{array}
\label{tilb}
\eeq
It follows from Proposition \ref{e*} and (\ref{tilb})
that the condition $\vep^*_i(\vec x)\leq \lan h_i,\lm\ran$
is equivalent to
$-\vp(\vec x)\leq \lan h_i,\lm\ran$ for any $\vp\in\Xi_{\io}^{(i)}$.
\qed

\newtheorem{cor4}[thm4]{Corollary}
\begin{cor4}
\label{cor-blm}
Let $\io$ be the same one as in Theorem \ref{thm-e*} and $\lm$
be a dominant integral weight. Then we have:
\beq
B(\lm)\cong {\rm Im}(\plm)=
\{\vec x\in \Sigma_{\io}\ot R_{\lm}\,|\,
\lan h_i,\lm\ran+ \vp(\vec x)\geq 0\,\,{\rm for \,\,any}\,\,i\in I\,\,{\rm and\,\,}\vp\in \Xi_{\io}^{(i)}\}.
\label{blm-e*}
\eeq
\end{cor4}

\vskip2pt
Furthermore, we also obtain the following combinatorial expression for the weight multiplicities
and the tensor-product multiplicities as follows:
The weight function of $\ZZ^{\ify}_{\io}[\lm]$ is 
described explicitly by (\ref{wt-vep-vp}):
$wt(\vec x)=\lm-\sum_{k}x_k\al_{i_k}$.
Set $W(\lm):=\{\nu\in P\,|\,B(\lm)_{\nu}\ne \emptyset\}$ and denote
the weight multiplicity of $\nu$ in $B(\lm)$ by $M_{\lm,\nu}$.
Any $\nu\in W(\lm)$ is in the form $\lm-\sum_{i}m_i\al_i$ $(m_i\in\ZZ_{\geq0})$. Then we have
\begin{cor4}
For $\nu=\lm-\sum_{i}m_i\al_i\in W(\lm)$, the weight multiplicity of $\nu$ is given by
\beq
M_{\lm,\nu}=\sharp\{\vec x\in \wtil B(\lm)\,|\,
m_i=\sum_{i_k=i}x_k\,\,{\rm for\,\,any\,\,}i\in I\}.
\label{wt-multi}
\eeq
\end{cor4}
Now, we describe so-called the Littlewood-Richardson number $c^{\nu}_{\lm,\mu}$.
More precisely, for dominant integral weights $\lm,\mu$ and $\nu$, let
$c^{\nu}_{\lm,\mu}$ be the number of irreducible components $V(\nu)$ in the
tensor product $V(\lm)\ot V(\mu)$.
Of course, it is same as the number of connected components $B(\nu)$
in tensor product $B(\lm)\ot B(\mu)$.
To do this we need the follwing proposition similar to Proposition 3.2.1\cite{N}:
\newtheorem{pro4}[thm4]{Proposition}
\begin{pro4}
\label{highest}
For  dominant integral weights $\lm$ and $\mu$, an element $u\ot v\in B(\lm)\ot B(\mu)$ satisfies
$\eit(u\ot v)=0$ for any $i\in I$ if and only if
$\eit u=0$ and $\eit^{\lan h_i,\lm\ran+1}v=0$ for any $i\in I$.
\end{pro4}

{\sl Proof.\,\,}
The argument in the proof of Proposition 3.2.1\cite{N} can be applied
to any integrable highest weight modules
for symmetrizable Kac-Moody Lie algebras. \qed

Here note that the condition $\eit^{\lan h_i,\lm\ran+1}v=0$ is equivalent
to the one $\vep_i(v)\leq\lan h_i,\lm\ran$ and the explicit form of $\vep_i$
is given in (\ref{wt-vep-vp}). Here we set
\beq
E^{(i)}:=\{\sigma_k(\vec x)\,:\,i_k=i\}\cup\{\sigma^{(i)}_0(\vec x)\}.
\label{EI}
\eeq
\begin{cor4}
For dominant integral weight $\lm,\mu$ and $\nu$,
we have
\beq
c^{\nu}_{\lm,\mu}=\sharp\{\vec x\in \wtil B(\mu)\,|\,
wt(\vec x)=\nu-\lm\,\,{\rm and\,\,}\zeta(\vec x)\leq \lan h_i,\lm\ran
\,\,{\rm for\,\,any\,\,}i\in I\,\,{\rm and}\,\,\zeta\in E^{(i)}\}.
\label{LR}
\eeq
\end{cor4}

\section{Rank 2 case}
\setcounter{equation}{0}
\renewcommand{\theequation}{\thesection.\arabic{equation}}

In this section, we apply  Theorem \ref{thm-e*} and 
Corollary \ref{cor-blm} to the case for the
Kac-Moody algebras of rank 2.
We adopt the same setting as in \cite[Sect.4]{NZ}.
Without loss of generality, we can and will assume that $I=\{1,2\}$,
and $\io = (\cd,2,1,2,1)$.
The Cartan data is given by:
$$
\lan h_1,\al_1\ran= \lan h_2,\al_2\ran=2, \,\, \lan h_1,\al_2\ran=-c_1,
\,\, \lan h_2,\al_1\ran=-c_2.
$$
Here we either have $c_1 = c_2 = 0$, or both $c_1$ and $c_2$ are
positive integers.
We set $X = c_1 c_2 - 2$, and define the integer sequence
$a_l = a_l (c_1, c_2)$ for $l \geq 0$ by setting $a_0 = 0, \, a_1 = 1$
and, for $k \geq 1$,
\beq
a_{2k}  = c_1 P_{k-1} (X), \,\,
a_{2k+1} = P_k (X) + P_{k-1} (X),
\label{defcoeff}
\eeq
where the $P_k (X)$ are {\it Chebyshev polynomials} given 
by the following generating function:
\beq
 \sum_{k \geq 0} P_k (X) z^k = (1 - X z + z^2)^{-1}.
\label{gen-cheb}
\eeq
Here define $a'_l(c_1,c_2):=a_l(c_2,c_1)$.
Let $l_{\rm max} = l_{\rm max} (c_1, c_2)$ be the minimal index
$l$ such that $a_{l+1} < 0$
(if $a_l \geq 0$ for all $l \geq 0$,
then we set $l_{\rm max} = + \infty$). 
By inspection, if $c_1 c_2 = 0$ (resp. $1,2,3$) 
then $l_{\rm max} = 2$ (resp. $3, 4, 6$).
Furthermore, if $c_1 c_2 \leq 3$ then $a_{l_{\rm max}} = 0$ and 
$a_l > 0$ for $1 \leq l < l_{\rm max}$.
On the other hand, if $c_1 c_2 \geq 4$, i.e., $X \geq 2$,
it is easy to see from (\ref{gen-cheb})
that $P_k (X) > 0$ for $k \geq 0$, hence 
$a_l > 0$ for $l \geq 1$; in particular, in this case
$l_{\rm max} = + \infty$.  

\newtheorem{thm6}{Theorem}[section]
\begin{thm6}
\label{rank 2-thm}
\begin{enumerate}
\item
In the rank 2 case, for a dominant integral weight 
$\lm=\lm_1\Lm_1+\lm_2\Lm_2$ $(\lm_1,\lm_2\in \ZZ_{\geq0})$ 
the image of the embedding $\Psi^{(\lm)}_{\io}$ 
is given by
\beq
{\rm Im} \,(\Psi^{(\lm)}_{\io}) = \left\{(\cd,x_2,x_1)\in\ZZ_{\geq0}^{\ify}: 
\begin{array}{l}
x_k = 0 \,\,{\rm for}\,\,k > l_{\rm max},\,\,\lm_1\geq x_1, \\
a_l x_l -a_{l-1} x_{l+1} \geq0 ,\\
\lm_2+a'_{l+1}x_l-a'_lx_{l+1}\geq0,\\
{\rm for}\,\,1 \leq l < l_{\rm max}
\end{array}
\right\}.
\label{rank 2-poly}
\eeq
\item
For any $b\in B(\ify)$, writing $\Psi_{\io}(b)=(\cd,x_2,x_1)$,
we have
\beq
\vep^*_1(b)=x_1,\q 
\vep^*_2(b)={\rm max}_{1\leq l\leq l_{\rm max}}
\{a'_lx_{l+1}-a'_{l+1}x_l\}.
\eeq
\end{enumerate}
\end{thm6}

\vskip5pt
{\sl Proof.}\,\,
In order to apply Corollary \ref{cor-blm},
we shall describe the set of linear functions
$\Xi_{\io}^{(\ify)}$ and $\Xi_{\io}^{(i)}$, and check 
that $\io$ satisfies the strict positivity assumption.


The set $\Xi_{\io}^{(\ify)}=\Xi_{\io}$ has been given in \cite[Lemma 4.2]{NZ}. 
In particular, it is shown that the positivity assumption are satisfied.
Hence, by Lemma 4.2 in \cite{NZ}, we have

\newtheorem{lm6}[thm6]{Lemma}
\begin{lm6}
\label{rank2-lm1}
\begin{enumerate}
\item [{\rm (i)}]
For $k \geq 1$ and $0 \leq l < l_{\rm max}$, we set
\beq
\vp^{(l)}_k = S_{k+l-1}\cd S_{k+1}S_{k} x_k;
\label{S-forms}
\eeq
in particular, $\vp^{(0)}_k = x_k$.  
 If $k$ is odd then 
$\vp^{(l)}_k = a_{l+1} x_{k+l} - a_l x_{k+l+1}$; 
if $k$ is even then 
$\vp^{(l)}_k = a'_{l+1} x_{k+l} - a'_l x_{k+l+1}$.
\item [{\rm (ii)}] If $c_1 c_2 \leq 3$, i.e., $l_{\rm max} < + \infty$,
then $\vp^{(l_{\rm max}-1)}_k  = -x_{k+l_{\rm max}}$.
\item [{\rm (iii)}] The set $\Xi^{(\ify)}_{\io}$ consists of all linear forms
$\vp^{(l)}_k$ with $k \geq 1$ and $0 \leq l < l_{\rm max}$.
\item [{\rm (iv)}] The positivity assumption for the sequence $\io$ is satisfied.
\end{enumerate}
\end{lm6}

\vskip5pt

Now we return to the proof of Theorem \ref{rank 2-thm}.
It is remained to describe
$$
\Xi^{(i)}_{\io}=\{S_{j_k}\cd S_{j_1}\xii(\vec x)
\,\,|\,\,k\geq0,\,j_1,\cd,j_k\geq1\}.
\q(i=1,2)
$$
Here we see the explicit form of $\xii$:
\beq
\xi^{(1)}=-x_1,\qq \xi^{(2)}=c_2x_1-x_2.
\label{xi12}
\eeq
It is evident that
$\Xi^{(1)}_{\io}=\{-x_1\}$.
The proof of the theorem is completed by the following lemma:
\begin{lm6}
\label{xi2}
\begin{enumerate}
\item [{\rm (i)}]
For $1 \leq l < l_{\rm max}$, we set
\beq
\eta_{l} = S_{l-1}\cd S_{2}S_{1}(\xi^{(2)});
\label{xi-forms}
\eeq
in particular, $\eta_{1} = \xi^{(2)}=c_2x_1-x_2$.  
Then we have $\eta_{l}=a'_{l+1}x_{l}-a'_{l}x_{l+1}$.
\item [{\rm (ii)}] If $c_1 c_2 \leq 3$, i.e., $l_{\rm max} < + \infty$,
then $\eta_{l_{\rm max}-1}  =-x_{l_{\rm max}}$.
\item [{\rm (iii)}] The set $\Xi^{(2)}_{\io}$ consists of all linear forms
$\eta_{l}$ with $1 \leq l < l_{\rm max}$.
\item [{\rm (iv)}] Any element in $\Xi^{(2)}_{\io}\setminus\{\xi^{(2)}\}$ has non-negative coefficients
for $x_1$ and $x_2$.
\end{enumerate}
\end{lm6}

{\sl Proof.\,}
We can check (ii) by direct calculations for $c_1c_2=0,1,2,3$. The statement
(iv) is immediate from (i) and (iii). Thus we shall show (i) and (iii).
Since $a'_l\geq0$, we have
\beqnn
S_{2k}(a'_{2k+1}x_{2k}-a'_{2k}x_{2k+1}) & = & 
a'_{2k+1}x_{2k}-a'_{2k}x_{2k+1}-a'_{l+1}\beta_{2k}\\
& = & (c_2a'_{2k+1}-a'_{2k})x_{2k+1}-a'_{2k+1}x_{2k+2}\\
& = & a'_{2k+2}x_{2k+1}-a'_{2k+1}x_{2k+2},
\eeqnn
where we use the relation $a'_{2k+2}=c_2a'_{2k+1}-a'_{2k}$.
Thus we get $S_{2k}\eta_{2k}=\eta_{2k+1}$. Similarly, we obtain 
$S_{2k-1}\eta_{2k-1}=\eta_{2k}$, 
$S_{2k+1}\eta_{2k}=\eta_{2k-1}$ and
$S_{2k}\eta_{2k-1}=\eta_{2k-2}$.
We also have $S_j\eta_{k}=\eta_{k}$ if $j\ne k,k+1$.
These imply (i) and also (iii).\qed

\vskip3pt
Applying Lemma \ref{rank2-lm1} and Lemma \ref{xi2} 
to Corollary \ref{cor-blm} 
we conclude that
\beq
{\rm Im} \,(\Psi^{(\lm)}_{\io}) =
\{(\cd,x_2,x_1)\in\ZZ^{\ify} \,|\,
\hspace{-5pt}
\begin{array}{l}
\vp^{(l-1)}_k(\vec x) \geq 0, \,\,
\lm_1\geq x_1\,\,{\rm and}\,\,\lm_2\geq -\eta_{l}(\vec x)\\
{\rm for}\,\, k \geq 1, \, 1 \leq l < l_{\rm max}
\end{array}\}
\label{intermed}
\eeq
Comparing (\ref{intermed}) with the desired answer (\ref{rank 2-poly}), and using 
parts (i) and (ii) of Lemma \ref{rank2-lm1},
it only remains to show that the inequalities $\vp^{(l)}_k \geq 0$
in (\ref{intermed}) are redundant when $k > 1$ and $l < l_{\rm max}-1$,
that is, they are consequences of the remaining inequalities.
This can be shown by the same way as in the proof of Theorem 4.1 in
\cite{NZ}.

The proof of (ii) is evident from Theorem \ref{thm-e*} and
Lemma \ref{xi2}.
\qed

%
\vskip5pt

Note that the cases when $l_{\rm max} < +\infty$, or equivalently, 
the image ${\rm Im} \,(\Psi_{\io})$ is contained in a lattice of finite rank,
just correspond to the Lie algebras 
$\ge=$ $A_1 \times A_1$, $A_2$, $B_2$ or $C_2$, $G_2$.

In conclusion of this section, we illustrate Theorem \ref{rank 2-thm}
by the example when $c_1 = c_2 = 2$, i.e., $\ge$ is the affine 
Lie algebra of type $A^{(1)}_1$, following to \cite{NZ}.
In this case, $X = c_1 c_2 - 2 = 2$.
It follows at once from (\ref{gen-cheb}) that $P_k (2) = k+1$;
hence, (\ref{defcoeff}) gives $a_l = l$ for $l \geq 0$.
We see that for type $A^{(1)}_1$, 
$$
B(\lm)\cong 
{\rm Im} \,(\Psi^{(\lm)}_{\io}) = \{(\cd,x_2,x_1)\in\ZZ_{\geq0}^{\ify}: 
\begin{array}{l}
l x_l - (l-1) x_{l+1} \geq0, \,\,\lm_1\geq x_1\,\,{\rm and}\,\,\\
\lm_2+(l+1)x_l-lx_{l+1}\geq0 \,\,{\rm for}\,\, l \geq 1
\end{array}
\}, 
$$
and for $\vec x=(\cd, x_2,x_1)\in \Sigma_{\io}$ we have
$$
\vep_1^*(\vec x)=x_1,\,\,\,
{\rm and}\,\,\,
\vep_2^*(\vec x)=
{\rm max}_{l\geq 1}
\{lx_{l+1}-(l+1)x_l\}.
$$

\section{$A_n$-case}
\setcounter{equation}{0}
\renewcommand{\theequation}{\thesection.\arabic{equation}}
We shall apply Theorem \ref{thm-e*} and 
Corollary \ref{cor-blm} to the case when $\ge$ is of type $A_n$. 
Let us identify the index set $I$ with $[1,n] := \{1,2,\cd,n\}$ 
in the standard way; thus, the Cartan matrix 
$(a_{i,j}= \lan h_i,\al_j\ran )_{1 \leq i,j \leq n}$ is given by 
$a_{i,i}=2$, $a_{i,j}=-1$ for $|i-j|=1$, and 
$a_{i,j}=0$ otherwise. 
As the infitite sequence $\io$ let us take 
the following periodic sequence 
$$
\io = \cd,\underbrace{n,\cd,2,1}_{},
\cd,\underbrace{n,\cd,2,1}_{},\underbrace{n,\cd,2,1}_{}.
$$

Following to \cite[Sect.5]{NZ}, we shall 
change the indexing set for $\ZZ^{\ify}$
from $\ZZ_{\geq 1}$ to $\ZZ_{\geq 1} \times [1,n]$, which is given by 
the bijection 
$\ZZ_{\geq 1} \times [1,n] \to \ZZ_{\geq 1}$ 
($(j;i) \mapsto (j-1)n + i$). 
According to this, we will write an element $\vec x \in \ZZ^{\ify}$
as a doubly-indexed family $(x_{j;i})_{j \geq 1, i \in [1,n]}$.
We will adopt the convention that $x_{j;i} = 0$ unless
$j \geq 1$ and  $i \in [1,n]$; in particular, $x_{j;0} = x_{j;n+1} = 0$
for all $j$. 

\newtheorem{thm7}{Theorem}[section]
\begin{thm7} 
\label{A_n}
Let $\lm=\sum_{1\leq i\leq n}\lm_i\Lm_i$ $(\lm_i\in \ZZ_{\geq0})$ 
be a dominant integral weight. 
\begin{enumerate}
\item
In the above notation, the image ${\rm Im} \,(\Psi^{(\lm)}_{\io})$ 
is the set of all integer families $(x_{j;i})$ such that 
\beqn
&&\hspace{-30pt}\hbox{
$x_{1;i} \geq x_{2;i-1} \geq \cd \geq x_{i;1} \geq 0$ 
for $1 \leq i \leq n$}
\label{sl-1}\\
&&\hspace{-30pt}\hbox{
$x_{j;i} = 0$ for $i+j > n+1$, }
\label{j;i=0}\\
&&\hspace{-30pt}\hbox{
$\lm_i\geq x_{j;i-j+1}-x_{j;i-j}$ for
$1\leq j\leq i\leq n$.}
\label{sl-2}
\eeqn
\item
For any $b\in B(\ify)$, writing $\Psi_{\io}(b)=(\cd,x_2,x_1)$ 
we have
\beq
\vep^*_i(b)={\rm max}_{1\leq j\leq i}
\{x_{j;i-j+1}-x_{j;i-j}\}.
\eeq
\end{enumerate}
\end{thm7}

\vskip5pt

{\sl Proof.}\,\,
We will follow the proof of Theorem \ref{rank 2-thm}.
So we first describe the set of linear functions
$\Xi_{\io}^{(i)}$ $(i=1,\cd,n,\ify)$, and check that $\io$ 
satisfies the strict positivity assumption.
As in the previous section, 
we set
\beqnn
\Xi^{(\ify)}_{\io} & := & 
\{S_{j_l}\cd S_{j_1}x_{j_0}\,\,|\,\,l\geq 0,j_0,j_1,\cd,j_l\geq1\},\\
\Xi^{(i)}_{\io} & := & 
\{S_{j_k}\cd S_{j_1}\xi^{(i)}(\vec x)\,\,|\,\,k\geq 0,j_0,j_1,\cd,j_l\geq1\}\,\,(i\in I),\\
\eeqnn

The explicit description of the 
set $\Xi_{\io}^{(\ify)}=\Xi_{\io}$ is given in \cite[Lemma 5.2]{NZ} and it is shown that 
the sequence $\io$  satisfies the positivity assumption, that is, in this setting
any linear form $\vp=\sum_k\vp_{j;i}x_{j;i}\in \Xi_{\io}$ has the property
$\vp_{1;i}\geq0$ for any $i=1,\cd,n$. 
Then (\ref{sl-1}) follows from Theorem 6.1 in \cite{NZ}.

Therefore, in order to complete the proof, 
it is sufficient for us to show that 
\beq
\Xi^{(i)}_{\io}=\{- x_{j;i-j+1}+x_{j;i-j}\,\,|\,\,1\leq j\leq i\},
 \,\,(i\in I).
\label{A_n-xi-i}
\eeq
Let us write $F^{(i)}$ for the R.H.S. of (\ref{A_n-xi-i}). 
By the definition, we have $\xi^{(i)}=-x_{1;i}+x_{1;i-1}$ and 
then $\xi^{(i)}\in F^{(i)}$. 
Here for $(j;i) \in \ZZ_{\geq 1} \times [1,n]$, we will write 
the piecewise-linear transformation $S_{(j-1)n+i}$ 
as $S_{j;i}$; if $(j;i) \notin \ZZ_{\geq 1} \times [1,n]$ 
then $S_{j;i}$ is understood as the identity transformation. 
By the direct calculations, we obtain immediately,
\beqn
&&\hspace{-35pt}S_{p;q}(- x_{j;i-j+1}+x_{j;i-j})\nn\\
&&\hspace{-25pt}=\cases{
-x_{j+1,i-j}+x_{j+1,i-j-1}& if $(p;q)=(j;i-j)$ and $j<i$,\cr
-x_{j-1,i-j+2}+x_{j-1,i-j+1}& if $(p;q)=(j;i-j+1)$ and $j\ne 1$,\cr
- x_{j;i-j+1}+x_{j;i-j} & otherwise,\cr}
\label{Spq-sl}
\eeqn
where note that if $j=i$, $- x_{j;i-j+1}+x_{j;i-j}=-x_{i;1}$. 
This implies that $F^{(i)}$ is closed by the action of $S_{p;q}$ and 
all elements  are obtained from $\xi^{(i)}$, which shows (\ref{A_n-xi-i}).
The strict positivity assumption follows from (\ref{A_n-xi-i}) 
immediately. Thus, by virtue of Theorem \ref{thm-e*} 
and Corollary \ref{cor-blm} if $\vec x=\Psi_{\io}(b)$ we have
$$
\vep^*_i(b)={\rm max}
\{ x_{j;i-j+1}-x_{j;i-j}\,\,|\,\,1\leq j\leq i\},
$$
which implies (ii) and then we have (i).\qed
%

\vskip 5pt

As we mentioned in 2.4, we give 
the example which doe not satisfy the positivity assumption.
\newtheorem{ex6}[thm7]{Example}
\begin{ex6}
We consider the case $\ge=\ssl_4$ and take 
the sequence $\io=\cd 2\,1\,2\,3\,2\,1$,
where we do not need the explicit form of ``$\cd$'' in $\io$.
For the simplicity, we write $\vec x=(\cd,x_2,x_1)$ for
an element $\vec x\in \ZZ^{\ify}_{\io}$.
In this setting, we have $\beta_1=x_1-x_2-x_4+x_5$, $\beta_2=x_2-x_3+x_4$ and 
$5^{(-)}=1$. Then 
$ S_1(x_1)=x_1-\beta_1=x_2+x_4-x_5,$ 
$S_2S_1(x_1)=x_2+x_4-x_5-\beta_2=x_3-x_5$ and 
$S_5S_2S_1(x_1)=x_3-x_5+\beta_1=x_1-x_2+x_3-x_4$.
Thus we see the form $S_5S_2S_1(x_1)$ has the negative coefficient for $x_2$,
which breaks the positivity assumption.
Furthermore, this case is not ample.
Fix $\lm\in P_+$ with $\lan h_2,\lm\ran >0$. 
Since $\beta^{(-)}_2=-\lan h_2,\lm\ran+x_2-x_1$ and 
$\what S_5\what S_2\what S_1(x_1)=S_5S_2S_1(x_1)$, 
we have
$$
\what S_2\what S_5\what S_2\what S_1(x_1)=x_1-x_2+x_3-x_4+\beta^{(-)}_2
=-\lan h_2,\lm\ran+x_3-x_4,
$$
which implies 
$\vec 0=(\cd,0,0)\not\in\Sigma_{\io}[\lm]$ 
because of $\lan h_2,\lm\ran>0$.
\end{ex6}

\section{$A^{(1)}_{n-1}$-case}
\setcounter{equation}{0}
\renewcommand{\theequation}{\thesection.\arabic{equation}}
In this section we shall treat 
the affine Lie algebra $\ge=A^{(1)}_{n-1}$.
We will assume that $n \geq 3$ since the case of $A^{(1)}_{1}$
was already treated in Sect.5.
As in \cite{NZ} we will identify the index set $I$ with $[1,n]$
in the way such that the Cartan matrix 
$(a_{i,j}= \lan h_i,\al_j\ran )_{1 \leq i,j \leq n}$ is given by 
$a_{i,i}=2$, $a_{i,j}=-1$ for $|i-j|=1$ or $|i-j|=n-1$, and 
$a_{i,j}=0$ otherwise. 
As the infinite sequence we take the following 
periodic sequence 
$$
\io = \cd,\underbrace{n,\cd,2,1}_{},
\cd,\underbrace{n,\cd,2,1}_{},\underbrace{n,\cd,2,1}_{}.
$$


In the rest of this section, we will use the notation (\cite{NZ}):
$$
j;i[k] := k-1 + (j-1)(n-1) + i.
$$
Thus, the correspondence $(j;i) \mapsto j;i[k]$ is a bijection from 
$\ZZ_{\geq 1} \times [1,n-1]$ to $\ZZ_{\geq k}$. 
If there is no confusion, we shall use $j;i$ for $j;i[1]$.
This bijection transforms the usual linear order on $\ZZ_{\geq k}$
into the {\it lexicographic} order on $\ZZ_{\geq 1} \times [1,n-1]$
given by
$$(j';i') < (j;i) \,\,{\rm if}\,\, j' < j \,\,{\rm or} \,\, 
j'=j,\, i' < i.$$
As in \cite[Sect 6]{NZ} we consider integer ``matrices"  $C = (c_{j;i})$ indexed by
$\ZZ_{\geq 1} \times [1,n-1]$, and such that $c_{j;i} = 0$
for $j \gg 0$.
With every such $C$ and any $k \geq 1$ we associate a linear form
$\vp_{C[k]}$ on $\ZZ^{\ify}$ given by
$\vp_{C[k]} = \sum_{j;i} c_{j;i} x_{j;i[k]}.$

\nd
For any $(j;i) \in \ZZ_{\geq 1} \times [1,n-1]$, we set
$s_{j;i} = s_{j;i}(C) = c_{1;i} + c_{2;i} + \cd + c_{j;i}.$

\newtheorem{df7}{Definition}[section]
\begin{df7}[\cite{NZ}]
A integer matrix $C$ indexed by $\ZZ_{\geq1}\times [1,n-1]$
(and each of the corresponding forms $\vp_{C[k]}$)
 is called {\it admissible} if it satisfies the following conditions 
(same as (6.2)--(6.5) in \cite{NZ}):
\beq
s_{j;i} \geq 0 \,\,{\rm for}\,\, (j;i) \in \ZZ_{\geq 1} \times [1,n-1].
\label{adm1}
\eeq
\beq
s_{j;i} = \delta_{i,1} \,\,{\rm for}\,\, j \gg 0. 
\label{adm2}
\eeq
\beq
\sum_{(j';i') \leq (j;i)} s_{j';i'} \leq j
\,\,{\rm for \,\, any}\,\, (j;i), 
\,\,{\rm with \,\, the \,\, equality \,\, for}
\,\, j \gg 0.
\label{adm3}
\eeq
\beq
{\rm If}\,\, s_{j;i} > 0 \,\,{\rm then}\,\, s_{j';i'} > 0
\,\,{\rm for \,\, some}\,\, (j';i') \,\,{\rm with}\,\, 
(j;i) < (j';i') \leq (j+1;i).
\label{adm4}
\eeq
\end{df7}

Let us denote the set of all admissible matrices by ${\cal C}$ and 
$C_0$ for the matrix given by $c_{j;i}=\del_{j;i,1;1}$.
Then we have $\vp_{C_0 [k]} = x_k$. The following lemma is shown in 
\cite{NZ}, which is used repeatedly in the subsequent arguments.
\newtheorem{lm7}[df7]{Lemma}
\begin{lm7}[Lemma 6.3 \cite{NZ}]
\label{c0}
The matrix $C_0$ is the only admissible matrix with $c_{1;1}=s_{1;1}>0$.
\end{lm7}

\newtheorem{thm8}[df7]{Theorem}
\begin{thm8}
\label{affine A}
For $\lm=\sum_i\lm_i\Lm_i\in P_+$ and the sequence $\io$ as above,
we have
\beq
{\rm Im}\,(\Psi^{(\lm)}_{\io})
=\left\{\vec x\in \ZZ^{\ify}[\lm]\,\left|\right.\,
\begin{array}{l}
\vp_{C[k]} (\vec x) \geq 0 {\rm \,\,for \,\,any\,\,} C
\in {\cal C}\,\,{\rm and\,\,} k\geq1,\\
\lm_i\geq x_{j;i}-x_{j;i-1}{\rm \,\,for \,\,} j
\geq 1{\rm \,\,and\,\,}1\leq i\leq n-1,\\
\lm_n+\vp_{C[0]}(\vec x)\geq0{\rm \,\,for \,\,any\,\,}C
\in {\cal C}\setminus \{C_0\}.
\end{array}
\right\}.
\label{thm-affine}
\eeq
\end{thm8}
Here note that we treat the matrix $C[0]$ 
in the third condition of (\ref{thm-affine}).
In this case there is no object corresponding to $\vp_{C_0[0]}=x_0$, 
but it is removed from ${\cal C}$. 
Furthermore, by Lemma \ref{c0} the matrix
with non-trivial $c_{1;1}$ is only $C_0$.
Thus the R.H.S of (\ref{thm-affine}) is well-defined.

\vskip5pt
{\sl Proof.\,\,}
Let $\Xi_k$ be the set of linear forms obtained 
by applying $S_j$'s on the linear form $x_k$ as in \cite[Sect 6]{NZ} and 
we denote $\Xi_k$ by $\Xi^{(\ify)}_k$. 
Then by Lemma 6.2 in \cite{NZ}, 
we have
\beq
\Xi^{(\ify)}_k=\{\vp_{C[k]}(\vec x)\,|\,C\in {\cal C}\}.
\label{xi-c}
\eeq
In order to complete the proof of the theorem, 
it suffices to show the following:
\newtheorem{pr7}[df7]{Proposition}
\begin{pr7}
\label{lm1 aff}
We have 
\beqn
\Xi^{(1)}_{\io}& = & \{-x_{1;1}\},
\label{xii-aff0}\\
\Xi^{(i)}_{\io}& = &\{-x_{j;i}+x_{j;i-1}\,|\,j\geq 1\},
\,\,\,(1<i\leq n-1),
\label{xii-aff1}\\
\Xi^{(n)}_{\io}&=&\{\vp_{C[0]}(\vec x)\,|\,C\in {\cal C}\setminus \{C_0\}\}.
\label{xii-aff2}
\eeqn
\end{pr7}

\vskip 5pt 
{\sl Proof of Proposition \ref{lm1 aff}}\,\,\,\,
The proof of (\ref{xii-aff0})  is trivial.

\nd
Let us show (\ref{xii-aff1}).
Write $F^{(i)}$  for the R.H.S. of (\ref{xii-aff1}) 
($1<i\leq n-1$).
Using the double index $j;i$, the linear form $\beta_{j;i}$ 
can be written explicitly in the following form:
\beq
\begin{array}{lll}
\beta_{j;i}(\vec x) & = & x_{j;i}-x_{j;i+1}-x_{j+1,i}+x_{j+1,i+1},\\
\beta_{(j;i)^{(-)}}(\vec x) & = & 
\cases{
x_{j-1;i-1}-x_{j-1;i}-x_{j,i-1}+x_{j,i} & if $(j;i)>(2;1)$,\cr
0 & if $(1;1)\leq (j;i)\leq (2;1)$.\cr}
\end{array}
\label{betaji}
\eeq
Here note that $x_{j;n}$ means $x_{j+1;1}$,  $x_{j;0}$ means
$x_{j-1;n-1}$ if $j>1$ and $x_{j,i}$ means 0 if $j\leq 0$.
which is the different convention from the $A_n$-case.
By the definition of $\xi^{(i)}$ we have
$\xi^{(i)}(\vec x)=-x_{1;i}+x_{1;i-1}$ ($1<i\leq n-1$).
Then $\xi^{(i)}(\vec x)\in F^{(i)}$. 
By using the explicit form of $\beta_{j;i}$ in (\ref{betaji}),
 we obtain the similar formula to (\ref{Spq-sl}):
\beq
S_{p;q}(- x_{j;i}+x_{j;i-1})=
\cases{
-x_{j+1,i}+x_{j+1,i-1}& if $(p;q)=(j;i-1)$,\cr
-x_{j-1,i}+x_{j-1,i-1}& if $(p;q)=(j;i)>(2;1)$,\cr
- x_{j;i}+x_{j;i-1} & otherwise.\cr}
\label{Spq-aff}
\eeq
This implies that any form in $F^{(i)}$ is generated from $\xi^{(i)}$ and 
the set $F^{(i)}$ is closed under the action of $S_{p;q}$.

\nd
Before showing (\ref{xii-aff2}) we see the following lemma:
\begin{lm7}
\label{c22}
Suppose that $C=(c_{j;i})\in {\cal C}\setminus\{C_0\}$ satisfies $c_{2;2}<0$.
Then we have $c_{1;2}=c_{2;1}=1$, $c_{2;2}=-1$ and $c_{j;i}=0$ for
other $j;i$.
\end{lm7}

{\sl Proof of Lemma \ref{c22}.\,\,\,}
By the definition of $s_{j;i}$, we have 
$c_{2;2}=s_{2;2}-s_{1;2}$. Thus, our assumption $c_{2;2}<0$ implies 
\beq
0\leq s_{2;2}< s_{1;2}.
\label{2;2-0}
\eeq
It is obtained by (\ref{adm3}) and (\ref{adm4}) that 
\beqn
&& s_{1;1}+s_{1;2}+\cd+ s_{1;n-1}\leq 1,
\label{2;2-1}\\
&& s_{1;1}+s_{1;2}+\cd+ s_{1;n-1}+s_{2;1}\leq 2,
\label{2;2-2}\\
&&\hbox{
$s_{j';i'}>0$ for some $(1;2)<(j';i')\leq (2;2)$}.
\label{2;2-3}
\eeqn
Since $C\ne C_0$, by Lemma \ref{c0}
we have $c_{1;1}=s_{1;1}=0$. 
Then by (\ref{adm1}), (\ref{2;2-0}) and (\ref{2;2-1}),
we get 
\beqn
&&s_{1;2}=c_{1;2}=1
\label{2;2-7}\\
&&s_{1;3}=s_{1;4}=\cd= s_{1;n-1}=s_{2;2}=0.
\label{2;2-4}
\eeqn
Then we have 
\beq 
c_{2;2}=s_{2;2}-s_{1;2}=-1.
\label{2;2-5}
\eeq
Furthermore, by (\ref{2;2-2}), (\ref{2;2-3}) and (\ref{2;2-4}) we have
\beq
s_{2;1}=1.
\label{2;2-6}
\eeq
Here we need the following lemma to complete the proof of Lemma \ref{c22}:
\begin{lm7}
\label{lm2;2}
If $s_{j;1}=1$ and $\sum_{(j';i')\leq (j;1)}s_{j';i'}=j$, we have
$s_{k;i}=\del_{i,1}$ for $(k;i)\geq (j;1)$.
\end{lm7}

{\sl Proof.\,\,}
It follows from (\ref{adm3}) and (\ref{adm4}) that 
\beqn
&&\hbox{
$\sum_{(j';i')\leq (j;n-1)}s_{j';i'}\leq j$,}
\label{2;2-8}\\
&&\hbox{
$\sum_{(j';i')\leq (j+1;1)}s_{j';i'}\leq j+1$,}
\label{2;2-9}\\
&&\hbox{
$s_{j';i'}>0$ for some $(j;1)<(j';i')\leq (j+1;1)$.}
\label{2;2-10}
\eeqn
Then by applying the assumption of the lemma to (\ref{2;2-8})--(\ref{2;2-10})
we have
$s_{j;2}=s_{j;3}=\cd=s_{j;n-1}=0$ and $s_{j+1;1}=1$,
which implies 
$\sum_{(j';i')\leq (j+1;1)}s_{j';i'}=j+1$. These are the assumption of the lemma
replaced $j$ by $j+1$. 
Therefore, the induction proceeds and then we get 
$s_{k;2}=s_{k;3}=\cd=s_{k;n-1}=0$ and $s_{k;1}=1$ for $k\geq j$.
\qed

\vskip5pt
By (\ref{2;2-2}), (\ref{2;2-7}) and (\ref{2;2-6}) we have 
$s_{2;1}=1$ and $\sum_{(j';i')\leq (2;1)}s_{j';i'}=2$. 
Then by using Lemma \ref{lm2;2}, we have
$s_{k;i}=\del_{i,1}$ for $(k;i)\geq (2;1)$ and then 
$$
c_{j;i}=s_{j;i}-s_{j-1;i}=\cases{
1 & if $(j;i)=(1;2)$, $(2;1)$,\cr
-1& if $(j;i)=(2;2)$,\cr
0 & otherwise,\cr}
$$
which is the desired result.
Then we finished the proof of Lemma \ref{c22}.
\qed

Let us show (\ref{xii-aff2}).
We write $F^{(n)}$ for the R.H.S. of (\ref{xii-aff2}). 
The explicit form of $\xi^{(n)}$ is 
\beq
\xi^{(n)}(\vec x)=x_{1;1}+x_{1;n-1}-x_{2;1}=x_{1;2[0]}+x_{2;1[0]}-x_{2;2[0]}.
\label{lmn-form}
\eeq
Observing this form carefully, formally we can write
\beq
\xi^{(n)}(\vec x)=S_{1;1[0]}x_{1;1[0]}.
\eeq
Of course there is nothing corresponding to $x_{{1;1}[0]}=x_0$. 
But, formally we have
\beq
F^{(n)} =  \{S_{j_l;i_l[0]}\cd S_{j_1;i_1[0]}(S_{1;1[0]}x_{1;1[0]})\,|\,
l\geq0,\,\,j_k\geq 1,\,\,i_k\in I\}
\label{lmn}
\eeq
By Lemma \ref{c0}, we know that the  form 
$\vp_{C[0]}(\vec x)=\sum c_{j;i}x_{j;i[0]}$
($C\in {\cal C}$) satisfying $c_{1;1}\ne0$
is only $x_{1;1[0]}=x_0$ corresponding to the matrix $C_0$.
Moreover, note that only the explcit form of 
$S_{2;2[0]}$ is different from 
those of $S_{2;2[k]}$ ($k\geq1$), namely, 
$(2;2[0])^-=(2;1)^-=0$ and $(2;2[k])^-=1;1[k]$ for $k\geq1$ and 
then $\beta_{2;2[0]^-}=0$ and $\beta_{2;2[k]^-}=\beta_{1;1[k]}$
for $(k\geq1)$.
But for $k\geq0$ 
the form $\vp_{C[k]}=\sum c_{j;i}x_{j;i[k]}$ with $c_{2;2}<0$ 
must be $x_{1;2[k]}+x_{2;1[k]}-x_{2;2[k]}$ by Lemma \ref{c22}.
Then we have that if $c_{2;2}<0$, 
\beq
S_{2;2[k]}\vp_{C[k]}(\vec x)=
\cases{\vp_{C[0]}=\xi^{(n)}=S_{1;1[0]}x_{1;1[0]}, & if $k=0$,\cr
       \vp_{C_0[k]}=x_k, & if $k\geq1$,\cr}
\label{220-22k}
\eeq
which implies that $\vp_{C_0[0]}=x_0$ never occurs in $F^{(n)}$ and 
$F^{(n)}$ is stable by the actions of $S_{j;i[0]}$.
Thus, by Lemma 6.2 in \cite{NZ} we obtain 
$$
F^{(n)}= \{\vp_{C[0]}(\vec x)\,|\,C\in {\cal C}\setminus \{C_0\}\}.
$$
Now, we completed the proof of Proposition \ref{lm1 aff}.\qed

\vskip5pt

{\sl Proof of Theorem \ref{affine A}.\,\,}
Let us see that the above proposition and lemma  imply our theorem.
In view of Corollary \ref{cor-blm} and Proposition \ref{lm1 aff}, 
it suffices to show that
the sequence $\io$ satisfies the strict positivity assumption.
We have already shown this for $\Xi^{(\ify)}_k$ in \cite{NZ}. 
Next, we see $\Xi^{(i)}_{\io}$ $(i=1,2,\cd,n)$.
Since $\xi^{(1)}=-x_{1;1}$ and $\xii=-x_{1;i}+x_{1;i-1}$ for 
$1<i\leq n-1$, among 
$\bigcup_{1\leq i\leq n-1}(\Xi_{\io}^{(i)}\setminus\{\xii\})$ the 
only linear form which has non-trivial coefficients
for $x_{1;i}$ ($1\leq i\leq n-1$) or $x_{2;1}$ is 
just $-x_{2;2}+x_{2;1}$, which has a positive coefficient for $x_{2;1}$.
The remaining case  is $i=n$.
In this case, the explicit form of $\xi^{(n)}$ is given by
$\xi^{(n)}(\vec x)=x_{1;1}+x_{1;n-1}-x_{2;1}.$
So it sufficies to show that any linear form 
$\vp(\vec x)=\sum_{j;i}c_{j;i}x_{{j;i}[0]}\in 
\Xi_{\io}^{(n)}\setminus\{\xi^{(n)}\}$
satisfies that $c_{1;2},\cd, c_{1;n-1},c_{2;1},c_{2;2}\geq0$.
Those for $c_{1;2},\cd,c_{1;n-1}$ are trivial from 
$s_{1;i}=c_{1;i}$ and (\ref{adm1}).
If we assume $c_{2;1}<0$, we have $c_{1;1}=s_{2;1}-c_{2;1}>0$.
By Lemma \ref{c0}, we get $c_{j;i}=\del_{j;i,1;1}$.
Now since $C\ne C_0$, we have $c_{2;1}\geq0$.
If we assume tha $c_{2;2}<0$, we have $\vp_{C[0]}=\xi^{(n)}$ by
Lemma \ref{c22}. Thus, we obtain the strict positivity assumption and
then completed the proof of the theorem.\qed

\end{document}